\documentclass{article}
\usepackage{amssymb,amsmath,
 theorem,euscript}

\newcounter{sec}

\newcounter{punct}[sec]

\def\punct{\refstepcounter{punct}{\arabic{sec}.\arabic{punct}.  }}

\def\COUNTERS{\addtocounter{sec}{1}
              \setcounter{punct}{0}
          \setcounter{equation}{0}
          \setcounter{theorem}{0}
          \setcounter{problem}{0}
            \setcounter{Apunct}{0}
          }

\newtheorem{theorem}{Theorem}[sec]
\newtheorem{proposition}[theorem]{Proposition}

\newtheorem{lemma}[theorem]{Lemma}

\newtheorem{corollary}[theorem]{Corollary}
\newtheorem{observation}[theorem]{Observation}

\newtheorem{conjecture}[theorem]{Conjecture}
\newtheorem{question}[theorem]{Question}

\def\COUNTERS{\addtocounter{sec}{1}
              \setcounter{punct}{0}
          \setcounter{equation}{0}
          \setcounter{theorem}{0}
          }

\def\GL{\mathrm  {GL}}
\def\U{\mathrm  U}
\def\OO{\mathrm  O}
\def\Sp{\mathrm  {Sp}}
\def\SO{\mathrm  {SO}}

\def\Gr{\mathrm{Gr}}

\def\OSp{\mathrm {OSp}}

\def\GD{\mathbf {GD}}

\def\Mat{\mathrm{Mat}}

\def\Pfaff{\mathrm {Pfaff}}

\def\B{\mathfrak B}

\def\phi{\varphi}
\def\epsilon{\varepsilon}
\def\kappa{\varkappa}

\def\bSp{\mathbf {Sp}}

\def\SF{\boldsymbol{\mathcal{S}} \mathbf{F}}

\def\le{\leqslant}
\def\ge{\geqslant}

\renewcommand{\Re}{\mathop{\rm Re}\nolimits}

\def\la{\langle}
\def\ra{\rangle}

\def\pia{\pi_\downarrow}

\def\super{\mathrm{s}}

\newcommand{\im}{\mathop{\rm im}\nolimits}
\newcommand{\indef}{\mathop{\rm indef}\nolimits}
\newcommand{\dom}{\mathop{\rm dom}\nolimits}

\newcommand{\Ber}{\mathop{\rm Ber}\nolimits}
\newcommand{\graph}{\mathop{\rm graph}\nolimits}
\newcommand{\nul}{\mathop{\mathsf {null}}\nolimits}
\newcommand{\Lagr}{\mathop{\rm Lagr}\nolimits}
\newcommand{\spin}{\mathop{\rm spin}\nolimits}

\def\bF{\mathbf F}

\def\cA{\mathcal A}

\def\cM{\mathcal M}

\def\cO{\mathcal O}

\def\cR{\mathcal R}
\def\cS{\boldsymbol{\mathcal S}}

\def\cV{\mathcal V}
\def\cW{\mathcal W}

\def\frB{\mathfrak B}
\def\frC{\mathfrak C}
\def\frD{\mathfrak D}

\def\frR{\mathfrak R}

\def\fra{\mathfrak a}

\def\frg{\mathfrak g}

\def\fro{\mathfrak o}
\def\frp{\mathfrak p}

\def\frs{\mathfrak s}

\def\fros{\mathfrak{s}}

\def\R {{\mathbb R }}
 \def\C {{\mathbb C }}
  \def\Z{{\mathbb Z}}

\def\Q{{\mathbb Q}}
\def\A{{\mathbb A}}

 \def\ov{\overline}
\def\wt{\widetilde}
\def\wh{\widehat}

\renewcommand\emptyset{\varnothing}

\def\arr{\rightrightarrows}

\def\ev{{\mathrm{even}}}
\def\od{{\mathrm{odd}}}
\def\even{{\mathrm{even}}}
\def\odd{{\mathrm{odd}}}

\def\q{\quad}

\def\F{\mathbf F}

\def\b{\mathfrak b}

\def\sm{\smallskip}

\begin{document}

\begin{center}

\bf\Large Gauss--Berezin integral operators,\\
 spinors
over orthosymplectic supergroups,
\\ and Lagrangian super-Grassmannians

\sc\large

\bigskip

Yuri A. Neretin%
\footnote{Supported by the grants FWF, projects  P19064, P31591,
the grant  NWO.047.017.015, and the grant JSPS-RFBR-07.01.91209}

\end{center}

\bigskip

{\small
We obtain explicit formulas for
the spinor representation $\rho$
of the real orthosymplectic supergroup $\OSp(2p|2q,\R)$
by integral 'Gauss--Berezin' operators.
Next, we extend $\rho$  to a complex domain and get  a representation of a larger semigroup, which is
a counterpart of Olshanski subsemigroups in semisimple Lie groups.
Further, we show that $\rho$ can be extended to an operator-valued function
on a certain domain in  the Lagrangian super-Grassmannian
(graphs of elements of the supergroup $\OSp(2p|2q,\C)$ are Lagrangian super-subspaces) and show
that this function is a 'representation' in the following sense:
we consider Lagrangian subspaces as linear relations and
composition of two Lagrangian relations in general position
 corresponds to a product of Gauss--Berezin operators%
\footnote{This paper is an  extended variant of preprint {\tt https://arxiv.org/abs/0707.0570v3}
and a strongly revised version of my earlier preprint
\cite{Ner-preprint} exploring
 other realization
of spinors.}.

}

\bigskip


\section{Introduction}

\COUNTERS

In the present paper, we consider algebras, superalgebras, functional spaces over complex numbers $\C$ and in few cases over real numbers $\R$.
The transposition of matrices is denoted by $A\mapsto A^t$.
The symbol $1_n$ denotes the unit matrix of size $n$.

\sm

{\bf \punct Orthosymplectic spinors.%
\label{ss:1.1}}
Let $x_1$, \dots, $x_p$ be real variables, $\xi_1$, \dots, $\xi_q$
be Grassmann variables,
$
\xi_i\xi_j=-\xi_j \xi_i$ for all $i$, $j$ (in particular, $\xi_i^2=0$).
We consider  differential operators
\begin{align}
&1,\quad x_k x_l,\quad x_k\frac{\partial}{\partial x_l},\quad \frac{\partial}{\partial x_k}
\frac{\partial}{\partial x_l},\quad \xi_m \xi_n, 
\quad \xi_m \frac{\partial}{\partial \xi_n},
\quad \frac{\partial}{\partial \xi_m} \frac{\partial}{\partial \xi_n},
\label{eq:p0}
\\
&\qquad x_k \xi_m,\quad x_k\frac{\partial}{\partial \xi_m},
\quad \xi_m \frac{\partial}{\partial x_l},\quad 
\frac{\partial}{\partial x_k} \frac{\partial}{\partial \xi_m},
\label{eq:p1}
\end{align}
where $1\le k,l\le p$, $1\le m, n\le q$,
acting in the space of polynomials in $x_1$, \dots, $x_p$,
 $\xi_1$, \dots, $\xi_q$. Denote by $\cR$ the space  of all {\it complex}
linear combinations of such operators.
We say, that a parity
of a nonzero monomial of degree 2 in $x_k$, $\frac{\partial}{\partial \xi_l}$,
$\xi_m$, $\frac{\partial}{\partial \xi_n}$  is $\ov 1$ if it contains either one $\xi_k$ or one 
$\frac{\partial}{\partial \xi_l}$. Otherwise parity is $\ov 0$.
So monomials \eqref{eq:p0} have parity $\ov 0$ and \eqref{eq:p1}
parity $\ov 1$.
We also say that a nonzero linear combination of monomials of
parity $\ov 0$ (resp. $\ov 1$) has parity $\ov 0$ (resp. $\ov 1$). 
Let $u$, $v\in \cR$ have parities $p(u)$, $p(v)$. We define 
the supercommutator of operators $u$, $v$ by 
\begin{equation}
[u,v]_s:= u v- (-1)^{p(u)p(v)} vu
\label{eq:supercommutator}
\end{equation}
and extend this operation to the whole $\cR$ by bilinearity.

It is  easy to see that the space $\cR$ is closed with respect to the 
supercommutator, and therefore we get a Lie superalgebra,
it is isomorphic to a direct sum of the orthosymplectic
Lie superalgebra $\mathfrak{osp}(2p|2q)$ and a trivial one-dimensional 
Lie algebra $\C$.

Let us recall a definition of $\mathfrak{osp}(2p|2q)$.
Denote the block $(p+p)\times(p+p)$-matrix $\begin{pmatrix}0&1_p\\-1_p&0 \end{pmatrix}$
by $J$ and the block $(q+q)\times(q+q)$-matrix $\begin{pmatrix}0&1_q\\1_q&0 \end{pmatrix}$
by $I$. The orthosymplectic Lie superalgebra $\mathfrak{osp}(2p|2q)$
consists of complex block $(2p+2q)\times (2p+2q)$-matrices
$\begin{pmatrix} A&B\\C&D\end{pmatrix}$
satisfying the condition%
\footnote{$J$ is a canonical form of a non-degenerate skew-symmetric matrix, $I$
is a canonical form of a nondegenerate complex symmetric matrix of an even order.
For our aims, $I$ is more convenient than the unit matrix, which is more
natural for  general theory.}
$$
\begin{pmatrix}J&0\\0&I \end{pmatrix}\begin{pmatrix} A&B\\C&D\end{pmatrix}
+\begin{pmatrix} A^t&C^t\\-B^t&D^t\end{pmatrix}\begin{pmatrix}J&0\\0&I \end{pmatrix}=0.
$$
We say that parity of matrices of the form $\begin{pmatrix}*&0\\0&* \end{pmatrix}$
is $\ov 0$,  parity of $\begin{pmatrix}0&*\\ *&0 \end{pmatrix}$ is $\ov 1$,
and define a supercommutator by formula \eqref{eq:supercommutator}.

If $p=0$, then we get the usual spinor representation of the orthogonal Lie algebra
$\mathfrak{o}(2q,\C)$ in the Grassmann algebra consisting of
'functions' in variables $\xi_1$, \dots, $\xi_q$. If $q=0$, then we get
a representation of symplectic Lie algebra $\mathfrak{sp}(2q)$ (symplectic spinors or oscillator representation).
Spinors and symplectic spinors are distinguished objects of representation theory.
 Orthosymplectic spinors  were considered in numerous
works, for instance, Berezin \cite{Ber-super} (with the construction mentioned above), 
Serov \cite{Serov} (where the spinor representations of the orthosymplectic supergroups
 $\OSp(2p|r)$ were  obtained), and 
 \cite{Bar1}, \cite{Bar2},  \cite{China2}, \cite{China}, \cite{Cou},
 \cite{Fra}, \cite{Fur2}, \cite{Khu}, \cite{LI}, \cite{Nish}.
 
 \sm
 
 {\sc Remark on notation.}
 In notation $\mathfrak{osp}(2p|2q)$ for the orthosymplectic Lie  superalgebra, I firstly write $2p$ corresponding to the 'human' (real or complex) variables
 and the symplectic Lie algebra $\mathfrak{sp}(2p)$;
  $2q$ corresponds to Grassmann variables and the orthogonal Lie algebra $\fro(2q)$. In literature, our  $\mathfrak{osp}(2p|2q)$ can be denoted
 by $\mathfrak{osp}(2q|2p)$ or $\mathfrak{spo}(2p|2q)$.
 \hfill $\boxtimes$
 
 \sm

 In this paper we write explicit
formula for representation of the corresponding global object, which is larger
than the real supergroup $\OSp(2p|2q)$. We do not assume that the reader
is familiar with super-mathematics and discuss orthosymplectic 
spinors as a topic of analysis and as a story about integral operators. The text is self-closed,  de facto we  use a minimal version 
of language of super-algebra and super-analysis%
\footnote{Lie superalgebras  can have a life of their own, without supergroups, supermanifolds, superintegration, etc., see \cite{Kac0}, \cite{BL}, Chapter 1. We
prefer to discuss supergroups --- Lie superalgebras are actually present only in Section
\ref{s:howe}} and  follow
DeWitt's \cite{DeW} way --- to consider linear spaces (modules) over Grassmann algebra with infinite number of generators%
\footnote{Different authors have different points of view to a formalization of 'super-analysis'. DeWitt's book  was an object a justified critisism in \cite{Sch}, \cite{Mol}. 
Our work is far from subleties of analysis on supermanifolds,  translation of our results to the more common functorial language is more-or-less automatical.}.

\sm

{\bf\punct Berezin formulas.} Apparently,
first elements  of a
strange
  analogy between the spinor representation of 
the orthogonal groups and the oscillator representation of symplectic groups%
\footnote{'{\it Spinor representation}' of orthogonal group is a common term
(the representation was discovered by \'Elie Cartan \cite{Car}, 1913), the term '{\it oscillator representation}' has several 
synonyms,
 namely,
the Weil representation, the Shale--Weil representation, the
Segal--Shale--Weil representation, the harmonic representation,
 the metaplectic representation, the symplectic spinors. The term 'oscillator representation' was proposed 
 by  Irving Segal (who was the first who described this representation \cite{Seg}, 1959).  For further references,
 see \cite{Ner-book}, \cite{Ner-lectures}.}
 were observed by K.~O.~Friedrichs in the early  1950s,
 see \cite{Fri}. He considered  spinors over symplectic groups
$\Sp(2n,\R)$ as a  kind of a self-obvious object (obtained by an application of the Stone-von Neumann theorem) and initiated a discussion about their extension to the case $n=\infty$,  
    see some historical comments in \cite{Ner-Ber}.

 In the beginning of  1960s, Feliks Berezin   obtained explicit formulas
\cite{Ber-DAN}, \cite{Ber-second}
 for both the representations.
We briefly recall his results. First of all, let us
realize the {\it real} symplectic group $\Sp(2n,\R)$
 as the group of
{\it complex} $(n+n)\times(n+n)$ matrices
\begin{equation}
g=
\begin{pmatrix}
\Phi&\Psi\\
\ov\Psi&\ov\Phi
\end{pmatrix}
\label{eq:g-sp}
\end{equation}
 satisfying  the condition
\begin{equation}
g
\begin{pmatrix}0&1\\-1&0 \end{pmatrix}
g^t =
\begin{pmatrix}0&1\\-1&0 \end{pmatrix}
.
\label{eq:g-sp-1}
\end{equation}
Similarly, we realize the {\it real} orthogonal group $\OO(2n,\R)$ as
the group of {\it complex} matrices
\begin{equation}
g=
\begin{pmatrix}
\Phi&\Psi\\
-\ov\Psi&\ov\Phi
\end{pmatrix}
\label{eq:g-o}
\end{equation}
 satisfying
\begin{equation}
g
\begin{pmatrix}0&1\\1&0 \end{pmatrix}
g^t =
\begin{pmatrix}0&1\\1&0 \end{pmatrix}
.
\label{eq:g-o-1}
\end{equation}

The oscillator (spinor) representation of $\Sp(2n,\R)$ (see \eqref{eq:g-sp}, \eqref{eq:g-sp-1})
is realized by the following
integral operators $W(\cdot)$ 
\begin{multline}
W
\begin{pmatrix}
\Phi&\Psi\\
\ov\Psi&\ov\Phi
\end{pmatrix}
f(z)
=\pm(\det \Phi)^{-1/2}
\times\\\times
 \int_{\C^n} \exp \Bigl\{\frac 12
\begin{pmatrix}
z& \ov u
\end{pmatrix}
\begin{pmatrix}
\ov\Psi\Phi^{-1}  &(\Phi^{t})^{-1}\\ \Phi^{-1}& -\Phi^{-1}\Psi
\end{pmatrix}
\begin{pmatrix}
z^t\\ \ov u^t
\end{pmatrix}\Bigr\}
\,f(u)e^{-|u|^2}du\,d\ov u \label{eq:ber-1}
\end{multline}
in the space of holomorphic functions on $\C^n$. Here
 the symbol
$^t$ denotes the transposition of matrices; $z=\begin{pmatrix}
z_1&\dots&z_n\end{pmatrix}$,
 $u=\begin{pmatrix}u_1&\dots&u_n\end{pmatrix}$
are row vectors,
$$
\begin{pmatrix}
z& \ov u
\end{pmatrix}:=\begin{pmatrix}
z_1&\dots&z_n& \ov u_1&\dots&\ov u_n \end{pmatrix}
$$
also is a row vector, and the expression in the curly brackets is a product of a row vector,
matrix, and a column vector (i.e., the whole expression is a scalar). We write $W(\cdot)$ in honor of A.Weil.

On the other hand, Berezin obtained formulas for the spinor
representation of the group $\OO(2n,\R)$ (we realize $\OO(2n,\R)$ by matrices  \eqref{eq:g-o}, \eqref{eq:g-o-1}). Exactly,
 he wrote 
'integral operators'  of the form
\begin{multline}
\mathrm{spin}\begin{pmatrix}
\Phi&\Psi\\
\ov\Psi&\ov\Phi
\end{pmatrix}
f(\xi) = \pm(\det \Phi)^{1/2}
\times\\\times
 \int \exp\Bigl\{\frac 12
\begin{pmatrix}
\xi& \ov \eta
\end{pmatrix}
\begin{pmatrix}
-\ov\Psi\Phi^{-1}  & (\Phi^{t})^{-1}\\ -\Phi^{-1}& -\Phi^{-1}\Psi
\end{pmatrix}
\begin{pmatrix}
\xi^t\\ \ov \eta^t
\end{pmatrix}\Bigr\}
\,f(\eta) e^{-\eta\ov\eta^t}
 d\eta \, d\ov\eta
 \label{eq:ber-2}
,
\end{multline}
here $\xi=\begin{pmatrix}\xi_1&\dots&\xi_n\end{pmatrix}$,
$\eta=\begin{pmatrix}\eta_1&\dots&\eta_n\end{pmatrix}$
are row-matrices; $\xi_j$, $\eta_j$, $\ov\eta_j$ are
 anti-commuting variables.%
\footnote{Actually, in  both the cases Berezin considered
$n=\infty$.}$^{,}$\footnote{Such a
formula
 makes sense only for an open
dense subset in $\SO(2n,\R)\subset\OO(2n,\R)$;
 this 
produces some difficulties below.
Spinors over infinite-dimensional symplectic and orthogonal groups
were independently introduced by Shale and Stinespring \cite{Sha}, \cite{SS}. However,
the  starting point of super-analysis and standpoint of the present work were Berezin's formulas \eqref{eq:ber-1}, \eqref{eq:ber-2}.}. 
The integral
in the right-hand
side is the Berezin integral, see Section
\ref{s:orthogonal}.

In fact, Berezin in his  book \cite{Ber-second}%
\footnote{On 'intellectual history' of this book, its
origins and influence, see \cite{Ner-Ber}.}
and the subsequent work \cite{Ber-grass} conjectured that there
is the analysis of
Grassmann variables parallel to the usual analysis (and these works 
contain  important elements of this analysis). The book
includes parallel exposition of the bosonic Fock space (analysis 
in infinite number of complex variables)
and the fermionic Fock space
(Grassmann analysis in infinite number of variables),
 this parallel seemed mysterious. However, the strangest elements of this
analogy were  formulas (\ref{eq:ber-1}) and (\ref{eq:ber-2}). This pushed him at the  end of 60s -- beginning of 70s
to the invention of the `super-analysis',
 which mixes even
(complex or real) variables and odd (Grassmann) variables,
see \cite{KNV}, Sect. 5.

Apparently, the first case of joining of classical and Grassmann analysis was the
paper by Berezin  and
G.~I.~Kats \cite{BK}, 1970, where formal Lie supergroups
corresponding to Lie superalgebras
were introduced. At least in 1965 Lie superalgebras  appeared in algebraic topology
(see Milnor, J.~Moore
\cite{MM}, this work was one of starting points for Berezin and Kats). Starting 1971-73 Lie superalgebras and supersymmetries were used 
in the  quantum field theory (e.g., \cite{GL}, \cite{VA}, \cite{VS}, \cite{SSt},  'super' originates from this literature).

\smallskip


{\bf\punct Aims of this paper. Orthosymplectic spinors.} We wish
to unite  formulas (\ref{eq:ber-1})--(\ref{eq:ber-2})
 and to
write explicitly
 the representation of the supergroup%
 \footnote{Thus, this paper returns to the initial point
of super-analysis. It seems strange that formula (\ref{eq:3}) was
 not written by Berezin himself.
  The author obtained it after  reading the  posthumous uncompleted
  book \cite{Ber-super} of Berezin; in my opinion, 
 it contains traces of attempts to do this
  (see also  \cite{BT}).
In the present paper, we use  tools that were unknown to
 Berezin.
A straightforward
 extension of \cite{Ber-second} leads to cumbersome calculations.
Certainly, these difficulties were surmountable. In any case, the present paper is a kind of a 'lost chapter' of the Berezin's  book \cite{Ber-second}
and my book \cite{Ner-book}.}
$\OSp(2p|2q,\R)$,
the operators of the representation have the form
\begin{multline}
T(g)f(z,\xi)=
\iint
\exp\left\{\frac 12
\begin{pmatrix}
z&\xi&\ov u&\ov\eta
\end {pmatrix}
\cdot\frR(g)\cdot
\begin{pmatrix}
z^t\\ \xi^t\\ \ov u^t\\ \ov\eta^t
\end {pmatrix}
\right\}
\times\\ \times
f(u,\eta)\,
e^{-u \ov u^t-\eta\ov\eta^t}\,du\,d\ov u\,d\eta\,d\ov\eta
\label{eq:3}
,
\end{multline}
where $g$ ranges in the  supergroup $\OSp(2p|2q,\R)$ and $\frR(g)$ is a
certain block matrix of size $(p+q+p+q)$ 
composed of elements a supercommutative algebra $\cA$
(we prefer to think that $\cA$ is the Grassmann algebra with infinite number
of generators $\fra_j$), see Section \ref{s:operators}.

\smallskip


{\bf \punct Aims of this paper. Gauss--Berezin
integral operators.}
 In \cite{NNO}, \cite{Ner-boson},
\cite{Ner-fermion} it was shown that spinor and oscillator representations
are actually  representations of categories. We obtain a
(non perfect) super-counterpart of these constructions.
Let us explain this in more detail.

First, let us consider the 'bosonic' case. Consider
a symmetric
$(n+n)\times (n+n)$-matrix,
$$
S=\begin{pmatrix} A&B\\ B^t&C \end{pmatrix}
.$$
Consider a {\it Gaussian integral operator}%
\footnote{We use the symbol $\frB$ in
honour of Berezin.}
\begin{multline}
\frB_+\begin{bmatrix} A&B\\ B^t&C \end{bmatrix}
f(z)=\\=
\int_{\C^n} \exp \Bigl\{\frac 12
\begin{pmatrix}
z& \ov u
\end{pmatrix}
\begin{pmatrix} A&B\\ B^t&C \end{pmatrix}
\begin{pmatrix}
z^t\\ \ov u^t
\end{pmatrix}\Bigr\}
\,f(u)e^{-|u|^2}du\,d\ov u
 \label{eq:ber-3}
.
\end{multline}
These operators are more general than (\ref{eq:ber-1}),
it can be shown (this was observed by G.~I.~Olshanski) that the symmetric matrices 
$\begin{pmatrix}
\ov\Psi\Phi^{-1}  &(\Phi^{t})^{-1}\\ \Phi^{-1}& -\Phi^{-1}\Psi
\end{pmatrix}$ 
in
 Berezin's formula (\ref{eq:ber-1}) 
 are unitary.

 It can be readily checked that bounded Gaussian operators form a semigroup,
 which
 includes the group $\Sp(2n,\R)$; its algebraic structure
 is described below in Section \ref{s:symplectic}.

\smallskip

On the other hand,
 we can introduce {\it Berezin operators} that are fermionic counterparts of
 Gaussian operators.
Namely, consider a skew-symmetric matrix
$\begin{pmatrix} A&B\\ -B^t&C \end{pmatrix}$
and the integral operator
\begin{multline}
\frB_-\begin{bmatrix} A&B\\ -B^t&C \end{bmatrix}
f(\xi)=
\\=
\int \exp \Bigl\{\frac 12
\begin{pmatrix}
\xi& \ov \eta
\end{pmatrix}
\begin{pmatrix} A&B\\- B^t&C \end{pmatrix}
\begin{pmatrix}
\xi^t\\ \ov \eta^t
\end{pmatrix}\Bigr\}
\,f(u)e^{-\eta\ov\eta}d\eta\,d\ov \eta
 \label{eq:ber-4}
.
\end{multline}
These operators are more general than
(\ref{eq:ber-2}); it can be shown that
in Berezin's formula (\ref{eq:ber-2}) the skew-symmetric
 matrix
 $\begin{pmatrix}
-\ov\Psi\Phi^{-1}  & (\Phi^{t})^{-1}\\ -\Phi^{-1}& -\Phi^{-1}\Psi
\end{pmatrix}$
is  contained in the pseudo-unitary group $\U(n,n)$.

For Berezin operators in  general position, the product
of Berezin operators has a kernel of the same form. But sometimes
this is not the case. It is possible to improve
the definition
(see our Section \ref{s:orthogonal}), and
to obtain a semigroup of Berezin operators.

In this paper we  introduce 'Gauss--Berezin
integral operators', which unify
 bosonic 'Gaussian operators' and fermionic 'Berezin operators'.
Our main result is a construction of a canonical bijection
between the set of all Gauss--Berezin operators and a certain
domain in Lagrangian super-Grassmannian;
also we propose
a geometric interpretation of products of Gauss--Berezin
operators.

Formula (\ref{eq:3}) for
superspinor representation of the supergroup
 $\OSp(2p|2q)$ is a byproduct
of the geometric construction.


\smallskip

{\bf\punct Aims of this paper. Possible applications.}
The spinor and oscillator representations
are important at least for the  following two reasons.

First, they are a basic tool 
in representation theory of infinite-dimensional groups,
see, e.g., \cite{Seg}, \cite{Olsh-GB}, \cite{Ner-book}, in particular, for classical groups,
 the group of diffeomorphisms of the circle, loop groups%
 \footnote{The parabolic induction, which is a main tool for 
 construction of representations of semisimple Lie groups does not
 work in these situations.}.

Second, the Howe duality for spinors is an important topic of classical
representation theory, see, for instance,
\cite{How},  \cite{KV}, \cite{Ada}.

Possible applications of our work in similar directions are discussed in Section \ref{s:howe}.

\smallskip

{\bf\punct Structure of this paper.}
I  tried to write a
self-contained paper, no preliminary knowledge of the
super-mathematics or  representation theory is
 assumed. But this implies a necessity of various
 preliminaries, which I provide. Also, we try to minimize the vocabulary of this text.

We start with an exposition of  spinor
and oscillator representations
in Sections \ref{s:orthogonal} and \ref{s:symplectic}. These sections contain also
a discussion of Gaussian integral operators
and Berezin (fermionic Gaussian) operators (a detailed exposition of these topics
is contained in the books \cite{Ner-book} and \cite{Ner-lectures}.

In Sections \ref{s:group} -- \ref{s:relations} we discuss
supergroups $\OSp(2p|2q)$, super-Grassmannians and super-linear relations.
I am trying to use  minimally necessary tools,
these sections do not contain an introduction  to super-science and its basic definitions.
For  generalities of super-mathematics, see \cite{Bern}, \cite{Ber-super},
\cite{Man}, \cite{Lei}, \cite{DeW},
 \cite{CF}, \cite{Var}, \cite{BL}.

Section \ref{s:integral} contains a  discussion
of super-analogue of the Gaussian integral, this is a simple
imitation of  the well-known formula
\begin{equation}
\int_{\R^n} \exp\Bigl\{- \frac 12 x Ax^t+bx^t\Bigr\}
\,dx=(2\pi)^{n/2} \det (A)^{-1/2}
\exp\Bigl\{- \frac 12 bA^{-1} b^t\Bigr\}
.
\label{eq:gaint}
\end{equation}

  Apparently, these calculations are
 written somewhere, but I do not know  references.

The Gauss--Berezin integral  operators  are introduced in Section
\ref{s:operators}.
Our main construction is the canonical one-to-one
correspondence
between super-linear relations and Gauss--Berezin operators,
 which
is obtained  in Section  \ref{s:correspondence}. This
immediately produces a representation
 of real supergroups $\OSp(2p|2q,\R)$.
For $g$ being in an 'open dense' subset in the supergroup,
the
operators of the representation
are of the form (\ref{eq:3}).

In the last section we discuss some 
 open problems.


\smallskip

{\bf Acknowledgements.} I am grateful to D.~V.~Alekseevski, \fbox{A.~L.~Onishchik}, and
A.~S.~Losev  for explanations of
super-algebra and super-analysis. The preliminary variant
\cite{Ner-preprint}
of this paper was  revised after a discussion with
D.~Westra, who proposed numerous suggestions
to improve the text.


\section{A survey of orthogonal
spinors.  Berezin \\ operators and  Lagrangian
linear relations}

\label{s:orthogonal}

\COUNTERS

This Section  is subdivided into 3 parts.

In Part A we develop
the standard formalism
of Grassmann algebras $\Lambda_n$
 and define Berezin operators $\Lambda_n\to\Lambda_m$, which
are in a some sense  morphisms of
 Grassmann algebras (but  not morphisms in the
category of algebras!).

In Part B we describe the 'geometrical' category $\GD$,
which is equivalent to the category of Berezin operators.
Morphisms of the category $\GD$ are certain Lagrangian subspaces.

In Part C we describe explicitly the correspondence
between Lagrangian subspaces and Berezin operators.

\smallskip

In some cases we present proofs or   explanations,
for a coherent treatment
see \cite{Ner-book}, Chapter 2.

\begin{center}
\bf A. Grassmann algebras and Berezin operators.
\end{center}


\nopagebreak

{\bf\punct Grassmann variables and Grassmann algebra.%
\label{ss:grassmann-algebra}}
We denote by $\xi_1$, \dots, $\xi_n$ 
Grassmann variables,
$$
\xi_i\xi_j=-\xi_j\xi_i
,
$$
in particular, $\xi_i^2=0$.
 Denote by $\Lambda_n$ the algebra of
polynomials {\it with complex coefficients} in these variables, evidently,
 $\dim \Lambda_n=2^n$.
The monomials
\begin{equation}
\xi_{j_1}\xi_{j_2}\dots\xi_{j_\alpha}, \qquad
\text{where $\alpha=0,1,\dots,n$
and $j_1<j_2<\dots<j_\alpha$},
\label{eq:standard-basis}
\end{equation}
form a basis of $\Lambda_n$.
Below  elements of Grassmann algebra
 are called {\it  functions.}

\smallskip


{\bf\punct Derivatives.} We define {\it left
differentiations} in $\xi_j$ as usual.
Exactly, if $f(\xi)$
does not depend on $\xi_j$, then
$$
\frac{\partial}{\partial \xi_j} f(\xi)=0, \qquad
\frac{\partial}{\partial \xi_j} \xi_jf(\xi)=f(\xi)
.
$$

Evidently,
$$
\frac{\partial}{\partial \xi_k}
\frac{\partial}{\partial \xi_l}=-
\frac{\partial}{\partial \xi_l}
\frac{\partial}{\partial \xi_k},\qquad
\left(\frac{\partial}{\partial \xi_k}\right)^2=0
.$$


{\bf\punct Exponentials.}
Let $f(\xi)$ be an {\it even} function,
 i.e., $f(-\xi)=f(\xi)$. We define its
exponential as usual,
$$
\exp\bigl\{f(\xi)\bigr\}:=\sum_{j=0}^\infty
\frac 1{j!} \, f(\xi)^j
.$$
Even functions $f$, $g$ commute, $fg=gf$, therefore
$$
\exp\{f(\xi)+g(\xi)\}=\exp\{f(\xi)\}\cdot \exp\{g(\xi)\}
.$$

\smallskip


{\bf\punct Berezin integral.
\label{ss:berezin-integral}} Let $\xi_1$, \dots, $\xi_n$ be
Grassmann variables. The {\it Berezin integral}
$$
\int f(\xi)\,d\xi=\int f(\xi_1,\dots,\xi_n)\,d\xi_1\dots d\xi_n
$$
is a linear functional on $\Lambda_n$ defined
by
$$
\int \xi_1\xi_{2}\dots \xi_n\,d\xi_1\dots d\xi_n=1,
$$
 integrals of all  other
monomials are zero.

The following formula for integration by parts holds
$$
\int f(\xi)\cdot \frac{\partial}{\partial \xi_k} \,
g(\xi)\,d\xi=- \int \frac{\partial f(-\xi)}{\partial \xi_k}
\cdot g(\xi)\, d\xi
.
$$


{\bf\punct Integrals with respect 'to odd Gaussian measure'.%
\label{ss:berezin-integral-2}} Let $\xi_1$, \dots, $\xi_q$ be
 as above. Let $\ov\xi_1$, \dots, $\ov\xi_q$ be another
collection of Grassmann variables,
$$
\ov \xi_k\ov\xi_l=-\ov\xi_l\ov\xi_k,\qquad
\xi_k\ov\xi_l=-\ov\xi_l\xi_k
.
$$

We need the following Gaussian expression
\begin{multline*}
\exp\bigl\{-\xi\ov\xi^t\bigr\}:=
\exp\{-\xi_1\ov \xi_1-\xi_2\ov\xi_2-\dots\}
=\\=
\exp\{\ov\xi_1 \xi_1\} \exp\{\ov\xi_2 \xi_2\}\dots=
(1+\ov\xi_1\xi_1) (1+\ov\xi_2\xi_2)
 (1+\ov\xi_3\xi_3)\dots
\end{multline*}
Denote
$$
d\ov\xi\, d\xi=d\ov\xi_1\, d\xi_1\, d\ov\xi_2\,d\xi_2\dots
$$
 Therefore,
\begin{equation}
\int\Bigl( \prod_{k=1}^m
(\ov\xi_{\alpha_k}\xi_{\alpha_k})\Bigr)
\, e^{-\xi\ov\xi^t}\,
d\ov\xi\,d\xi
=1,
\end{equation}
and the integral is zero for all  other monomials.
For instance,
\begin{gather*}
\int \xi_1\xi_2\ov\xi_3 \, e^{-\xi\ov\xi^t}\,d\ov\xi\,d\xi=0,
\qquad  \int \, e^{-\xi\ov\xi^t}\, d\ov\xi\,d\xi=+1,
\\
\int \xi_1 \ov\xi_1 \ov\xi_{33}\xi_{33}
 \, e^{-\xi\ov\xi^t}\,
 d\ov\xi\,d\xi= - \int (\ov\xi_1 \xi_1)\,
 (\ov \xi_{33} \xi_{33})
\, e^{-\xi\ov\xi^t}\,
 d\ov\xi\,d\xi=
 -1.
\end{gather*}

 Evidently,
\begin{align*}
\int f(\ov\xi)\,\cdot\,\frac{\partial}{\partial \xi_k}
g(\xi) \, e^{-\xi\ov\xi^t}\, d\ov\xi\,d\xi
=
\int \ov\xi_k
 f(-\ov\xi) \,\cdot\,
g(\xi)\, e^{-\xi\ov\xi^t}\,d\ov\xi\,d\xi
,
\\
\int \frac{\partial}{\partial\ov\xi_k}
 f(\ov\xi)\,\cdot\,g(\xi)
\, e^{-\xi\ov\xi^t}\,
d\ov\xi\, d\xi
=\int f(-\ov\xi)\cdot\xi_k g(\xi)\,
 e^{-\xi\ov\xi^t}\,d\ov \xi\,d\xi
.
\end{align*}


{\bf\punct Integral operators.%
\label{ss:integral-operators-1}}
Now, consider Grassmann algebras $\Lambda_p$
and $\Lambda_{q}$ consisting of polynomials
in Grassmann variables $\xi_1$, \dots, $\xi_p$
and $\eta_1$, \dots, $\eta_q$, respectively. For a function
$K(\xi,\ov\eta)$ we define an {\it integral operator}
$$
A_K:\Lambda_q\to\Lambda_p
$$
by
$$
A_K f(\xi)=\int K(\xi,\ov\eta)f(\eta)\,
e^{-\eta\ov\eta^t}\,d\ov\eta\,d\eta
.
$$

\begin{proposition}
\label{pr:grassmann-operator-symbol}
 The map $K\mapsto A_K$ is a
 one-to-one correspondence of the set of all polynomials
$K(\xi,\ov\eta)$ and the set of all linear maps
$\Lambda_q\to\Lambda_p$.
\end{proposition}

\smallskip

A proof of Proposition \ref{pr:grassmann-operator-symbol} is 
trivial. Indeed, expand
$$
K(\xi,\eta)
=\sum a_{i_1,\dots,i_l\, j_1,\dots,j_l}
\xi_{i_1}\dots\xi_{i_k}
\ov\eta_{j_1}\dots\ov\eta_{j_l}
.
$$
Then $a_{\dots}$ are the matrix elements%
\footnote{up to  signs.}
of $A_K$
 in the standard basis (\ref{eq:standard-basis}).

\begin{proposition}
\label{pr:product-operators-grassmann}
If $A:\Lambda_q\to\Lambda_p$ is determined by
the kernel $K(\xi,\eta)$ and $B:\Lambda_p\to\Lambda_r$
is determined by the kernel $L(\zeta,\ov\xi)$,
then the kernel of $BA$ is
$$
M(\zeta,\ov\eta)=\int L(\zeta,\ov\xi) K(\xi,\ov\eta)
e^{-\xi\ov\xi^t}\,d\ov\xi\,d\xi
.
$$
\end{proposition}

\smallskip


{\bf\punct Berezin operators in the narrow sense.%
\label{ss:berezin-narrow}} A {\it Berezin operator
$\Lambda_q\to\Lambda_p$ in the narrow sense} is an
operator of the form
\begin{multline}
\B\begin{bmatrix}
A&B\\-B^t&C
\end{bmatrix}
f(\xi):=\\=
\int \exp\Bigl\{\frac12 \begin{pmatrix}\xi&\ov \eta \end{pmatrix}
\begin{pmatrix} A&B\\ -B^t& C\end{pmatrix}
\begin{pmatrix}\xi^t\\ \ov\eta^t \end{pmatrix}\Bigr\}
f(\eta) \,e^{-\eta\ov\eta^t}\,d\ov\eta\,d\eta
\label{eq:berezin-narrow}
,
\end{multline}
where $A=-A^t$, $C=-C^t$. Let us explain the notation.

\smallskip

1. $\begin{pmatrix}\xi&\eta \end{pmatrix}$ denotes the
row-matrix
$$
 \begin{pmatrix}\xi&\ov\eta \end{pmatrix}
:=\begin{pmatrix}\xi_1&\dots&\xi_p &\ov\eta_1&\dots&\ov\eta_q
 \end{pmatrix}
.
$$
 Respectively, $\begin{pmatrix}\xi^t\\
\ov\eta^t \end{pmatrix}$ denotes the transposed
column-matrix.

\smallskip

2. The $(p+q)\times(p+q)$-matrix
$\begin{pmatrix} A&B\\ -B^t& C\end{pmatrix}$ is
skew-symmetric. The whole expression for the kernel has the form
$$
 \exp\Bigl\{ \frac12 \sum_{k\le p, l\le p} a_{kl}\xi_k\xi_l
+\sum_{k\le p,m\le q} b_{km} \xi_k\ov\eta_m +
\frac 12 \sum_{m\le q, j\le q} c_{mj}\ov\eta_m\ov\eta_j
\Bigr\}
.$$

\smallskip

{\bf\punct Product formula.%
\label{ss:product-formula-grassmann}}

\begin{theorem}
\label{th:product-berezin}
Let
$$
\B[S_1]=\B\begin{bmatrix} P&Q\\ -Q^t& R\end{bmatrix}:\Lambda_p\to\Lambda_q
,\qquad
\B[S_2]=
\B\begin{bmatrix} K&L\\ -L^t& M\end{bmatrix}:\Lambda_q\to\Lambda_r
$$
be Berezin operators in the narrow sense. Assume $\det(1-MP)\ne 0$. Then
\begin{equation}
\B[S_2]\,\B[S_1]=\Pfaff
\begin{pmatrix} M&1\\-1&P\end{pmatrix}
\B[S_2\circ S_1]
\label{eq:product1}
,
\end{equation}
where
\begin{multline}
\label{eq:strange-berezin}
\begin{pmatrix} K&L\\ -L^t& M\end{pmatrix}
\circ \begin{pmatrix} P&Q\\ -Q^t& R\end{pmatrix}=
\\=
\begin{pmatrix}
K+LP(1-MP)^{-1}L^t & L(1-PM)^{-1}Q\\
-Q^t(1-MP)^{-1} L^t& R-Q^t(1-MQ)^{-1}MQ
\end{pmatrix}
.
\end{multline}
\end{theorem}

The symbol $\Pfaff(\cdot)$ denotes the Pfaffian, see
the next subsection.

\smallskip

A calculation is not difficult, see Subsection
\ref{ss:grassmann-gauss}.  The formula \eqref{eq:strange-berezin} is clarified below, Theorem \ref{th:GD-berezin-equivalence}.

\smallskip


{\bf\punct Pfaffians and odd Gaussian integrals.%
\label{ss:pfaffian}} Let $R$ be a skew-symmetric $2n\times 2n$
matrix. Its {\it Pfaffian} $\Pfaff(R)$ is defined by the
condition
$$
\frac 1{n!}\Bigl(\frac12\sum_{kl} r_{kl} \xi_k\xi_l\Bigr)^n
=\Pfaff(R)\,\xi_1\xi_2\dots\xi_{2n-1}\xi_{2n}
$$
(where $\xi_1$, \dots, $\xi_{2n}$ are Grassmann variables).
In other words,
$$
\Pfaff(R)=\frac 1 {n!} \int \Bigl(\frac 12 \sum_{kl} r_{kl} \xi_k\xi_l\Bigr)^n\,d\xi
=\int\exp\Bigl\{\frac 12 \xi R\xi^t\Bigr\}\,d\xi
.
$$
Recall that
$$
\Pfaff(R)^2=\det R.
$$

More generally, let $\xi_1$, \dots, $\xi_{2n}$, $\theta_1$, \dots, $\theta_{2n}$
be pairwise anticommuting variables. Then (see \cite{Ner-book}, Theorem II.4.4)
\begin{equation}
\int\exp\Bigl\{\frac 12 \xi R\xi^t+\sum_{j=1}^{2n} \theta_j\xi_j \Bigr\}\,d\xi
=\Pfaff(R) \exp\Bigl\{\frac12 \theta R^{-1}\theta^t  \Bigr\}.
\label{eq:odd-gauss}
\end{equation}
Theorem \ref{th:product-berezin} follows from this formula in a straightforward way.

\smallskip


{\bf\punct The definition of Berezin operators.%
\label{ss:berezin-general}} Theorem \ref{th:product-berezin}
suggests an extension of the definition of Berezin operators.
Indeed, this theorem is perfect for operators in general position.
If $\det(1-MP)=0$, then we get
 an indeterminate 
of the form $0\cdot\infty$ in the product formula
(\ref{eq:product1}), (\ref{eq:strange-berezin}).

For this reason, consider the cone%
\footnote{Below, a {\it cone} in $\C^N$
 is a subset invariant with respect
 to homotheties
$z\mapsto s z$, where $s\in\C$.}
$\mathcal{C}$ of
all the operators of the form $s\cdot \B\begin{bmatrix}A&B^t\\-B^t&C \end{bmatrix}$,
 where $s$ ranges in
$\C$. Certainly, the cone $\mathcal{C}$ is not closed.
Indeed,
\begin{align*}
&\lim_{\epsilon\to 0}
\epsilon \exp\Bigl\{ \frac 1\epsilon \xi_1\xi_2 \Bigr\}
=\xi_1\xi_2;
\\
&\lim_{\epsilon\to 0}
\epsilon^2 \exp\Bigl\{
\frac 1\epsilon (\xi_1\xi_2+\xi_2\xi_4) \Bigr\}
=\xi_1\xi_2 \xi_3\xi_4
.
\end{align*}
This suggests the following definition.

\smallskip

{\it Berezin operators} are operators
 $\Lambda_q\to\Lambda_p$, whose
kernels have the form (cf. \cite{SMD})
\begin{equation}
s\cdot
\prod_{j=1}^m (\xi u^t_j+\ov\eta v^t_j)\cdot
 \exp\Bigl\{\frac12 \begin{pmatrix}\xi&\ov \eta \end{pmatrix}
\begin{pmatrix} A&B\\ -B^t& C\end{pmatrix}
\begin{pmatrix}\xi^t\\ \ov\eta^t \end{pmatrix}\Bigr\}
\label{eq:ber-wide}
,
\end{equation}
where

\smallskip

1. $s\in\C$;

\smallskip

2. $u_j\in\C^p$, $v_j\in\C^q$ are row-matrices;

\smallskip

3. $m$ ranges in the set $\{0,1,\dots,p+q\}$;

\smallskip

4. $\begin{pmatrix} A&B\\ -B^t& C\end{pmatrix}$
is a skew-symmetric $(p+q)\times(p+q)$-matrix.


{\bf\punct  The space of Berezin operators.}\nopagebreak

\begin{proposition}
\label{pr:cone-closed} {\rm a)} The cone of Berezin operators is closed
in the space of all linear operators $\Lambda_q\to\Lambda_p$.

\smallskip

{\rm b)} Denote by $\Ber^{[m]}$ the set of all  Berezin operators
with a given number $m$ of linear factors in {\rm(\ref{eq:ber-wide})}.
Then the closure of $\Ber^{[m]}$ is
$$
\Ber^{[m]}\cup \Ber^{[m+2]}\cup \Ber^{[m+4]}\cup \dots
$$
\end{proposition}

 Therefore, the cone of all  Berezin operators
$\Lambda_q\to\Lambda_p$ consists of two components,
namely
$$
\bigcup\limits_{\text{$m$ is even}}
\Ber^{[m]}
\qquad\text{and}\qquad
\bigcup\limits_{\text{$m$ is odd}}
\Ber^{[m]}.
$$
They are closures of $\Ber^{[0]}$ and $\Ber^{[1]}$, 
respectively.
Also,   kernels $K(\xi,\ov\eta)$ of Berezin operators satisfies
$$
K(-\xi,-\ov\eta)=K(\xi,\ov\eta),
\qquad \text{or}\qquad
 K(-\xi,-\ov\eta)=-K(\xi,\ov\eta),
$$
respectively.

\smallskip

{\sc Remark.} Kernels $K(\xi,\eta)$ of Berezin operators $\Lambda_q\to
\Lambda_p$ are canonically defined polynomials in $\xi$, $\eta$. However,  expressions
(\ref{eq:ber-wide}) for the kernels of   Berezin operators are
non-canonical if  numbers $m$ of linear factors are $\ge 1$. For
instance,
$$
(\xi_1+\xi_{33})(\xi_1+7\xi_{33})=
6\xi_1\xi_{33},
\qquad
\xi_1\exp\{\xi_1\xi_2+\xi_3\xi_4\}
=\xi_1 \exp\{\xi_3\xi_4\}
.
$$
A quadratic form in the exponential in (\ref{eq:ber-wide}) also is not defined canonically.
\qquad $\boxtimes$


\smallskip

{\bf\punct Examples of Berezin operators.}

\smallskip

 a) The identity operator is a Berezin operator. Its kernel
is
$
\exp\bigl\{\sum \xi_j\ov\eta_j\}
$.

b) More generally,  operators with kernels $ \exp\bigl\{\sum_{ij}
b_{ij} \xi_i\ov\eta_j\bigr\}$ are   operators
$\Lambda_q\to\Lambda_p$ defined by the substitution
$
\eta_i=
\sum_j b_{ji}\xi_j .
$

\smallskip

c) The operator with kernel $\xi_1\exp\{\sum \xi_j\ov\eta_j\}$
is the operator $f\mapsto \xi_1f$.

\smallskip

e) The operator with kernel $\ov\eta_1\exp\{\sum \xi_j\ov\eta_j\}$
is the operator $f\mapsto \frac{\partial}{\partial \xi_1}f$.

\smallskip

f) The operator in $\Lambda_n$ determined by the kernel
$(\xi_1+\ov\eta_1)(\xi_2+\ov\eta_2)\dots$ is  Hodge
$\boldsymbol{\star}$-operator. Namely, let $i_1<i_2<\dots<i_k$ be a subset in $\{1,2, \dots,n\}$. Let 
$j_1<j_2<\dots<j_{n-k}$ be the complementary subset. Then
$$
\boldsymbol{\star}\bigl(\xi_1\xi_2\dots\xi_{i_k}\bigr)=\xi_{j_1}\dots\xi_{j_{n-k}}.
$$

\smallskip

g) The operator  with kernel $K(\xi,\eta)=1$
is the projection to the vector $f(\xi)=1\in \Lambda_n$.

\smallskip

h) The operator $\frB$ with the kernel
$$
\exp\Bigl\{\frac 12 a_{ij}\xi_i\xi_j + \sum \xi_j\ov\eta_j\Bigr\}
$$
is the multiplication operator
$$
\frB f(\xi)= \exp\Bigl\{\frac 12 a_{ij}\xi_i\xi_j\Bigr\}\, f(\xi)
.
$$

i) A product of Berezin operators is a Berezin operator (see below).

\smallskip

j) Operators of the spinor representation of $\OO(2n,\C)$ are
Berezin operators.


\smallskip

{\bf\punct Another definition of Berezin operators.%
\label{ss:berezin-operators-2}}
Denote by $\frD[\xi_j]$
the following operators in $\Lambda_p$
\begin{equation}
\frD[\xi_j]f(\xi):=
\Bigl(\xi_j+\frac\partial{\partial \xi_j}
\Bigr)\,f(\xi)
\label{eq:frD}
.
\end{equation}
If $g$, $h$ do not depend on $\xi_j$, then
$$
\frD[\xi_j]\bigl(g(\xi)+\xi_j h(\xi)\bigr)=
\xi_j g(\xi)+h(\xi)
$$

In the same way, we define the operators
$\frD[\eta_j]:\Lambda_q\to\Lambda_q$.
Obviously
\begin{equation}
\frD[\xi_i]^2=1,\qquad
\frD[\xi_i]\,\frD[\xi_j]=-
\frD[\xi_j]\,\frD[\xi_i],\quad\text{for $i\ne j$}
\label{eq:frDfrD}
.
\end{equation}

A {\it Berezin operator} $\Lambda_q\to\Lambda_p$
is any operator that can be represented in the form
\begin{equation}
\frD[\xi_{k_1}]\dots \frD[\xi_{k_\alpha}]
\cdot \frB\cdot \frD[\eta_{m_1}]\dots \frD[\eta_{m_\beta}]
,
\label{eq:ber-wide2}
\end{equation}
where

\smallskip

--- $\frB$ is a Berezin operator in the narrow sense,

\smallskip

--- $\alpha$ ranges in the set $\{0,1,\dots,p\}$
and $\beta$ ranges in $\{0,1,\dots,q\}$

\smallskip

By (\ref{eq:frDfrD}) we can assume $k_1<\dots<k_\alpha$,
$m_1<\dots< m_\beta$.

\begin{proposition}
\label{pr:two-definitions}
 The two definitions of
Berezin operators are equivalent.
\end{proposition}

{\sc Remark.} Usually, a Berezin operator admits
 many representations in 
the form (\ref{eq:ber-wide2}).  In fact, the space of all
Berezin operators is a smooth cone, and
formula (\ref{eq:ber-wide2}) determines $2^{p+q}$ coordinate
systems on this cone. Any collection of $2^{p+q}-1$ charts have 
a non-empty complement.
\hfill $\boxtimes$

\smallskip


{\bf\punct The category of Berezin operators.
Groups of automorphisms.%
\label{ss:berezin-category}}

\begin{theorem}
\label{th:product-berezin-2}
Let
$$
\frB_1:\Lambda_q\to\Lambda_p,\qquad
 \frB_2:\Lambda_p\to\Lambda_r
 $$
be Berezin operators. Then $\frB_2\frB_1:\Lambda_q\to\Lambda_r$ is a
Berezin operator. If $\frB$ is an invertible Berezin operator, then
$\frB^{-1}$ is a Berezin operator.
\end{theorem}

By Theorem \ref{th:product-berezin-2} we get a category,
whose objects are Grassmann algebras $\Lambda_0$, $\Lambda_1$,
$\Lambda_2$, \dots and whose  morphisms are Berezin
operators.

Denote by $G_n$ the group of all  invertible Berezin operators
$\Lambda_n\to\Lambda_n$. By the definition, it contains the group
$\C^\times$ of all  scalar non-zero operators.

\begin{theorem}
\label{th:automorphisms}
$G_n/\C^\times\simeq\OO(2n,\C)/\{\pm 1\}$.
\end{theorem}

Here $\OO(2n,\C)$ denotes the usual group
of orthogonal transformations in $\C^{2n}$.

The map sending each element of   $\OO(2n,\C)/\{\pm 1\}$.
Moreover, this isomorphism sending an orthogonal matrix to the corresponding Berezin operator in $\Lambda_n$ (determined up to a factor) is nothing
but the spinor representation of
$\OO(2n,\C)$.

Our next aim is to describe explicitly the
category of Berezin operators. In fact, we intend to clarify
the strange matrix multiplication (\ref{eq:strange-berezin}).

\smallskip

\begin{center}
\bf B. Linear relations and the category $\GD$.
\end{center}

\smallskip


{\bf\punct Linear relations.%
\label{ss:linear-relations}} Let $V$, $W$ be linear spaces over
$\C$. A {\it linear relation} $P:V\arr W$ is a linear subspace
$P\subset V\oplus W$.

\smallskip

{\sc Remark}. Let $A:V\to W$ be a linear operator. Its {\it graph}%
\index{graph}
\label{N:graph}
 $\graph(A)\subset V\oplus W$ consists
of all  vectors $v\oplus Av$.
By the definition, $\graph(A)$ is a linear relation,
$$\dim \graph(A)=\dim V.$$


{\bf\punct Product of linear relations.%
\label{ss:prl}} Let $P:V\arr W$, $Q:W\arr
Y$ be linear relations. Informally, a product $QP$ of linear
relations is a product of many-valued maps. If $P$ takes a vector
$v$ to a vector $w$ and  $Q$ takes the vector $w$ to a vector $y$,
then $QP$ takes $v$ to $y$.

Now, we present a formal definition.
The {\it product} $QP$ is a linear relation $QP:V\arr Y$
consisting of all $v\oplus y\in V\oplus Y$ such that there exists
$w\in W$ satisfying  $v\oplus w\in P$, $w\oplus y\in Q$.

\smallskip

In fact, the multiplication of linear relations  extends the
usual matrix multiplication.

\smallskip


{\bf\punct Imitation of some standard definitions of
matrix theory.%
\label{ss:imitation}}

\smallskip

$1^\circ$. The {\it kernel} $\ker P$ consists of all $v\in V$
such
that $v\oplus 0\in P$. In other words,
$$\ker P=P\cap (V\oplus 0).$$

\smallskip

$2^\circ$. The {\it image} $\im P\subset W$ is the projection of
$P$ to $0\oplus W$.

\smallskip

$3^\circ$.  The {\it domain}
 $\dom P\subset V$ of
$P$ is the projection of $P$ to $V\oplus 0$.

\smallskip

$4^\circ$.  The {\it indefinity} $\indef P\subset W$ of $P$ is
$P\cap (0\oplus W)$.

\smallskip

The  definitions of a kernel and an image extend the
corresponding definitions for linear operators.
The definition of a domain
extends the usual definition of the domain of an unbounded
operator in an infinite-dimensional space.
For a linear operator, $\indef=0$.

\smallskip


{\bf\punct Lagrangian Grassmannian and orthogonal
groups.%
\label{ss:lo}} First, let us  recall some definitions.
Let $V$ be a linear space equipped with a non-degenerate symmetric
(or skew-symmetric) bilinear form $M$. A subspace $H$ is
{\it isotropic} with respect to the bilinear form $M$,
if $M(h,h')=0$ for all $h$, $h'\in H$. The  dimension of
an isotropic subspace satisfies
$$
\dim H\le \frac 12 \dim V
.
$$
A {\it Lagrangian subspace}%
\footnote{Mostly, the term 'Lagrangian subspace' is used for
spaces equipped with a skew-symmetric bilinear forms;
however, my usage is a usual slang.}
 is an isotropic subspace
whose dimension is precisely $\frac 12\dim V$.
By $\Lagr(V)$ we denote the {\it Lagrangian Grassmannian}, i.e.,
the space of all  Lagrangian subspaces in
$V$.

Consider a space $\cV_{2n}=\C^{2n}$ equipped with the
symmetric bilinear form $L=L_n$ determined by the matrix
$\begin{pmatrix}0&1_n\\1_n&0\end{pmatrix}$. Denote by $\OO(2n,\C)$ the
group of all  linear transformations $g$ of $\C^{2n}$
preserving $L$, i.e., $g$ must satisfy the condition
$$
g\begin{pmatrix}0&1_n\\1_n&0\end{pmatrix}g^t
=
\begin{pmatrix}0&1_n\\1_n&0\end{pmatrix}
.$$

Let $m$, $n=0$, 1, 2, \dots. Equip the space
$\cV_{2n}\oplus\cV_{2m}$ with  the symmetric bilinear
form $L^\ominus$ given by
\begin{equation}
\!\!
L^\ominus(v\oplus w,v'\oplus w')=
L_n(v,v')-L_m(w,w')\,\,
,
\label{eq:ominus}
\end{equation}
where $v$, $v'\in\cV_{2n}$ and $w$, $w'\in\cV_{2m}$.

\begin{observation}%
\footnote{I use term 'Observation' for  statements, which are important for understanding and became trivial or semi-trivial being formulated.}
\label{obs:graph-orthogonal} Let $g$ be an operator
in $\C^{2n}$.
Then $g\in\OO(2n,\C)$ if and only if its graph
is an $L^\ominus$-Lagrangian
subspace  in $\cV_{2n}\oplus\cV_{2n}$.
\end{observation}

This is obvious.

\smallskip


{\bf\punct Imitation of orthogonal groups.
Category $\GD$.%
\label{ss:GD}}
Now we define the category $\GD$. The objects are
the spaces $\cV_{2n}$, where $n=0,1,2,\dots$.
There are two types of morphisms
$\cV_{2n}\to\cV_{2m}$:

\smallskip

a) $L^\ominus$-Lagrangian subspaces $P\subset
\cV_{2n}\oplus\cV_{2m}$;
we regard them as linear relations.

\smallskip

b) a distinguished morphism%
\footnote{It is not identified with any linear relation.} denoted by
 $\nul_{2n,2m}$.

\smallskip

Now we define a product of morphisms.

\smallskip

--- Product of $\nul$ and any morphism is $\nul$.

\smallskip

--- Let $P:\cV_{2n}\arr\cV_{2m}$, $Q:\cV_{2m}\arr \cV_{2k}$ be
Lagrangian linear
relations.

$\phantom{\text{---}}$  Assume that
\begin{equation}
\ker Q\cap\indef P=0\qquad \text{or, equivalently,\,\, $\im P+\dom Q=\cV_{2m}$.}
\label{eq:null-condition}
\end{equation}

$\phantom{\text{---}}$ Then $QP$ is the product of linear relations.

\smallskip

--- If the  condition (\ref{eq:null-condition})
is not satisfied, then $QP=\nul$.

\smallskip

\begin{theorem}
\label{ss:GD-correct}
The definition is self-consistent, i.e.,
a product of morphisms is a morphism
and the multiplication is associative.
\end{theorem}

At first glance the appearance of $\nul$ seems strange;  its
necessity will be transparent immediately ($\nul$ correspond to
zero operators in the next theorem).
 Also, $\nul$ will be the source of some our
difficulties below.

\begin{theorem}
\label{th:GD-berezin-equivalence} The category of Berezin
operators defined up to scalar factors and the category $\GD$ are
equivalent.
\end{theorem}

In fact, there is a map that takes each Lagrangian
linear relation
$P:\cV_{2n}\arr \cV_{2m}$ to a nonzero Berezin operator
$\spin(P):\Lambda_n\to\Lambda_m$ such that for each
$R:\cV_{2q}\arr\cV_{2p}$, $Q:\cV_{2p}\arr\cV_{2r}$
$$
\spin(Q)\spin(R)=\lambda(Q,R)\spin(QR),
$$
where $\lambda(Q,R)$ is a constant. Moreover,
$$
\text{$\lambda(Q,R)=0$ if and only if $QR=\nul$.}
$$

We describe the correspondence between the category of Berezin
operators and the category $\GD$ explicitly
in Theorems \ref{th:spinor-correspondence},
\ref{th:spinor-correspondence-2}.
 First, we need  some auxiliary facts
 concerning  Lagrangian
Grassmannians (for a  detailed introduction
to Lagrangian Grassmannians,
see \cite{Arn}, Sect. 43,  \cite{Ner-lectures}, Chapter 3).

\begin{center}
\bf C. Explicit correspondence.
\end{center}


{\bf\punct Coordinates on the Lagrangian Grassmannian.%
\label{ss:coo-lagrange0}}
Recall that $\cV_{2n}$ is the space $\C^{2n}=\C^n\oplus \C^n$ equipped with the bilinear form $L=L_n$ with the matrix $\begin{pmatrix}0&1_n\\1_n&0\end{pmatrix}$.
We write $\cV_{2n}=\C^{2n}$
as
$$
\cV_{2n}=
\cV_n^+\oplus \cV_n^-=\C^n\oplus\C^n
,$$
in this decomposition the summands are Lagrangian subspaces.

\begin{lemma}
\label{l:potapov} Let $H\subset \cV_{2n}$ be an $n$-dimensional
subspace such that 
$H\cap\cV_n^-=0$. Under this condition $H$ is the
graph of an operator 
$$T_H:\cV_n^+\to \cV_n^-.$$
The following
conditions are equivalent

\smallskip

--- the matrix $T_H$ is skew-symmetric;

\smallskip

--- the subspace $H$ is Lagrangian
\end{lemma}

{\sc Proof.}
We write the bilinear form $L$ as
$$
L(v,w)=
L(v^+\oplus v^-,w^+\oplus w^-)=
v^+(w^-)^t +  w^+(v^-)^t
.
$$
By the definition, $v\in  H$ if and only if $v^-= v^+ T_H$.
For $v$, $w\in P$, we evaluate
$$
L(v^+\oplus v^+T,w^+\oplus w^+T)
=
v^+ ( w^+ T_H)^t+ v^+ T_H (w^+)^t
=v^+ (T_H+T_H^t) w^+
.
$$
Now the statement becomes obvious.
\hfill $\square$

\smallskip

Lemma \ref{l:potapov} defines a coordinate system in
$\Lagr(\cV_{2n})$. Certainly, this coordinate system does not
cover the whole space $\Lagr(\cV_{2n}) $.

\smallskip


{\bf\punct Atlas on the Lagrangian Grassmannian.%
\label{ss:atlas-lagrange0}}
Denote by $e_1^+$, \dots, $e_n^+$,  $e_1^-$, \dots, $e_n^-$ the standard basis in $\C^{2n}=\C^n\oplus \C^n$.
Let $J$ be a subset in $\{1,2,\dots,n\}$. Denote by $\ov J$ its
complement. We define the subspaces
\begin{align*}
\cV_n^+[J]:=\bigl(\oplus_{j\in J} \C e_j^+\bigr)
\oplus
\bigl(\oplus_{j\notin J} \C e_j^-\bigr)
,
\\
\cV_n^-[J]:=\bigl(\oplus_{j\notin J} \C e_j^+\bigr)
\oplus
\bigl(\oplus_{j\in J} \C e_j^-\bigr)
,
\end{align*}
then $\cV_{2n}=\cV_n^+[J]\oplus\cV_n^-[J]$. Denote by $\cM[J]$ the
set of all  Lagrangian subspaces $H$ in $\cV_{2n}$ such that
$H\cap \cV_n^-[J]=0$. Such subspaces are precisely graphs of
symmetric operators $\cV_n^+[J]\to\cV_n^-[J]$.

\begin{proposition}
\label{pr:maps-orthogonal}
The $2^n$ charts $\cM[J]$ cover the whole Lagrangian Grassmannian
$\Lagr(\cV_{2n})$.
\end{proposition}


{\bf\punct Atlas on the Lagrangian Grassmannian.
 Elementary reflections.%
\label{ss:elementary}}
 We can describe the same maps in  a slightly
different way. For $i=1,2,\dots,n$, define an elementary reflection
$\sigma_i:\cV_{2n}\to \cV_{2n}$ by
\begin{equation}
\sigma_i e_i^+=e_i^-,\qquad
\sigma_i e_i^-=e_i^+,\qquad
\sigma_i e_j^\pm=e_j^\pm\qquad\text{for $i\ne j$}
\label{eq:reflection}
.
\end{equation}

\begin{observation}
\begin{equation}
\cM[J]=
\Bigl(\prod_{i\in J} \sigma_i\Bigr) \cM[\varnothing]
\label{eq:sigma-karta}
.\end{equation}
\end{observation}

\smallskip


{\bf\punct  Components of Lagrangian Grassmannian.%
\label{ss:components-orthogonal}}

\begin{observation}
The Lagrangian Grassmannian in the space $\cV_{2n}$
consists of two connected components.%
\footnote{Recall that $\cV_{2n}$ is equipped with a
symmetric bilinear form. The usual Lagrangian Grassmannian
discussed in the next section is connected.
The orthosymplectic Lagrangian Grassmannian
(see Section \ref{s:grassmannian})
 consists of two components.}
\end{observation}

We propose two proofs to convince the reader.

\smallskip

1. The group $\OO(n,\C)$ is dense in $\Lagr(\cV_{2n})$.
This group has two components.

\smallskip

2. It can be readily checked that
$\Lagr(\cV_{2n})$ is a homogeneous space
$$
\Lagr(\cV_{2n})\simeq\OO(2n,\C)/\GL(n,\C)
.
$$
 The group $\OO(2n,\C)$
consists of two components and the group $\GL(n,\C)$
is connected. \hfill $\square$

\smallskip

{\sc Remark.} Two components of the Lagrangian
Grassmannian correspond to two components of the
space of Berezin operators. \hfill $\boxtimes$

\smallskip


{\bf\punct Coordinates on the set of morphisms of $\GD$.%
\label{ss:coordinates-morphisms}}
We can apply the  reasoning of
Section \ref{ss:coo-lagrange0}
to 
$\Lagr(\cV_{2n}\oplus
\cV_{2m})$. Due to the minus in  formula
(\ref{eq:ominus}) we must  take care of the signs in
Lemma \ref{l:potapov}.

\begin{lemma}
\label{l:potapov-2}
Decompose $\cV_{2n}\oplus\cV_{2m}$
as%
\footnote{We emphasize that the source space $\cV_{2n}$ and the
target space $\cV_{2m}$ are mixed in the next row.}
$$
\cV_{2n}\oplus\cV_{2m}=
(\cV_n^-\oplus \cV_m^+)
\oplus
(\cV_n^+\oplus \cV_m^-)
.$$
Let $P$ be an $(m+n)$-dimensional subspace such that
$P\cap (\cV_n^+\oplus \cV_m^-)=0$,
i.e., $P$ is the graph of an operator
$$
\cV_n^-\oplus \cV_m^+
\to
\cV_n^+\oplus \cV_m^-
.$$
Then the following conditions are
equivalent

\smallskip

--- $P\in\Lagr(\cV_{2n}\oplus\cV_{2m})$;

\smallskip

--- $P$ is the graph of an operator
 having the form
\begin{equation}
\begin{pmatrix}A&B\\B^t& C\end{pmatrix},\qquad
\text{where $A=-A^t$, $C=-C^t$}
.
\label{eq:ABBC}
\end{equation}
\end{lemma}


{\bf\punct Creation-annihilation operators.%
\label{ss:creation-grassmann}}
Let $\cV_{2n}$ be  as above. Decompose
$$
\cV_{2n}=\cV_n^+\oplus \cV_n^-,
\qquad\text{where $\cV_n^+:=\C^n\oplus 0$,
 \, $\cV_n^-:=0\oplus\C^n$.}
$$
Let us write elements of $\cV_{2n}$ as
$$
v:=\begin{pmatrix} v_1^+&\dots&v_n^+&v_1^-&\dots&v_n^-\end{pmatrix}
.
$$
For each $v\in \cV_{2n}$, we define a creation-annihilation
operator $\wh a(v)$ in $\Lambda_n$ by
$$
\wh a(v)f(\xi):=\Bigl( \sum_j  v_j^+ \xi_j +\sum_j v_j^-\frac{\partial}{\partial \xi_j}
\Bigr) f(\xi)
.
$$
Evidently,
$$
\wh a(v) \wh a(w)+\wh a(w)\wh a(v)=L(v,w)\cdot 1
.$$


{\bf\punct A construction of the correspondence.%
\label{ss:functor-orthogonal}}
 Let $\B:\Lambda_q\to\Lambda_p$ be a
nonzero Berezin operator.
Consider the subspace $P=P[\B]\subset
\cV_{2q}\oplus \cV_{2p}$ consisting of $v\oplus w$ such that
\begin{equation}
\wh a(w)\,\B=\B \,\wh a(v)
\label{eq:defining-fermi}
.
\end{equation}

\begin{theorem}
\label{th:spinor-correspondence}
{\rm a)} $P[\B]$ is a morphism of the category
$\GD$, i.e., a Lagrangian subspace.

\smallskip

{\rm b)} The map $\B\mapsto P$ is a bijection
\begin{multline*}
\left\{
\begin{matrix}
\text{Set of nonzero Berezin operators $\Lambda_q\to\Lambda_p$}\\
\text{defined up to a scalar factor}
\end{matrix}\right\}
\longleftrightarrow
\\
\longleftrightarrow
\left\{
\begin{matrix}
\text{The set of non-'$\nul$' morphisms $\cV_{2q}\to\cV_{2p}$}\\
\text{of the category $\GD$}
\end{matrix}
\right\}.
\end{multline*}

\smallskip

{\rm c)} Let $\B_1:\Lambda_q\to\Lambda_p$, $\B_2:\Lambda_p\to\Lambda_r$
be Berezin operators. Then
\begin{equation}
\B_2\B_1=0
\quad
\text{if and only if}\quad P[\B_2]\,P[\B_1]=\nul
.
\label{eq:null}
\end{equation}

d) Otherwise,
$$
P[\B_2\B_1]=P[\B_2] P[\B_1]
.
$$
\end{theorem}

{\sc Sketch of proof.} 1.  First, let us consider a Berezin
operator in the narrow sense. We write the equation
(\ref{eq:defining-fermi})
\begin{multline*}
\!
\Bigl(\sum w_j^+ \xi_j+\sum w_j^-
\frac{\partial}{\partial \xi_j}\Bigr)\!
\int\! \exp\Bigl\{\frac12 \begin{pmatrix}\xi&\ov \eta \end{pmatrix}
\!
\begin{pmatrix} A&B\\ -B^t& C\end{pmatrix}
\!
\begin{pmatrix}\xi^t\\ \ov\eta^t \end{pmatrix}\Bigr\}
f(\eta) e^{-\eta\ov\eta^t}  d\ov \eta\,d\eta=
\\=
\int \exp\Bigl\{\frac12 \begin{pmatrix}\xi&\ov \eta \end{pmatrix}
\begin{pmatrix} A&B\\ -B^t& C\end{pmatrix}
\begin{pmatrix}\xi^t\\ \ov\eta^t \end{pmatrix}\Bigr\}
\Bigl(\sum v_k^+ \eta_k+\sum v_k^-
\frac{\partial}{\partial \eta_k}\Bigr)
f(\eta) e^{-\eta\ov\eta^t} \, d\ov \eta\,d\eta
,
\end{multline*}
or, equivalently,
\begin{multline*}
\Bigl(\sum w_j^+ \xi_j+\sum w_j^-
\frac{\partial}{\partial \xi_j}
-\sum v_k^- \ov\eta_k-\sum v_k^+ \frac{\partial}{\partial \eta_k}
\Bigr)
\times\\ \times
\exp\Bigl\{\frac12 \begin{pmatrix}\xi&\ov \eta \end{pmatrix}
\begin{pmatrix} A&B\\ -B^t& C\end{pmatrix}
\begin{pmatrix}\xi^t\\ \ov\eta^t \end{pmatrix}\Bigr\}
=0
\end{multline*}
We differentiate the exponential
and get
\begin{multline*}
\Bigl[\sum_j w^+_j\xi_j+\sum_j w_j^-
\Bigl(\sum_i a_{ji}\xi_i+\sum_mb_{jk}\ov\eta_m\Bigr)
-
\sum_k v_k^- \ov\eta_k
-\sum_k v_k^+
\\+
\Bigl(-\sum_l b_{lk}\xi_l +\sum_m c_{km} \ov\eta_m  \Bigr)
\Bigr]
\times\\ \times
\exp\Bigl\{\frac12 \begin{pmatrix}\xi&\ov \eta \end{pmatrix}
\begin{pmatrix} A&B\\ -B^t& C\end{pmatrix}
\begin{pmatrix}\xi^t\\ \ov\eta^t \end{pmatrix}\Bigr\}
=0
.
\end{multline*}

Therefore,
$$
\begin{pmatrix}
v^+& w^-
\end{pmatrix}
=
\begin{pmatrix}
v^-& w^+
\end{pmatrix}
\begin{pmatrix} A&-B\\ -B^t& -C\end{pmatrix}
.
$$

By Lemma \ref{l:potapov-2}, $P$ is a  Lagrangian subspace in
$\cV_{2q}\oplus\cV_{2p}$.

\smallskip

2. Next, let $\frC$ be a Berezin operator having the form
(\ref{eq:ber-wide2}),
\begin{equation}
\frC=
\frD[\xi_{k_1}]\dots \frD[\xi_{k_\alpha}]
\cdot \B\cdot \frD[\eta_{m_1}]\dots \frD[\eta_{m_\beta}].
\label{eq:ber-wide3}
\end{equation}
We note that
$$
\wt a(v)\, \frD[\xi_i]=
\frD[\xi_i]\, \wt a(\sigma_i v)
,
$$
where $\sigma_i$ is an elementary  reflection
(\ref{eq:reflection}). Let the operator $\B$ satisfy
$$
\wh a(w)\,\B=\B\,\wh a(v)
,
$$
where $v\oplus w$ ranges in $P[\B]$. Then the operator $\frC$
satisfies
$$
\wh a(\sigma_{k_1}\dots\sigma_{k_\alpha} w)\,\frC =
\frC\,\wh a(\sigma_{m_1}\dots\sigma_{m_\beta} v)
.
$$
Therefore, the corresponding linear relation is
\begin{equation}
\sigma_{k_1}\dots\sigma_{k_\alpha}  P\sigma_{m_1}\dots\sigma_{m_\beta}
\label{eq:sigmaPsigma}
\end{equation}
and  is Lagrangian. This proves a).

\smallskip

3. By Proposition \ref{pr:maps-orthogonal}, the sets
(\ref{eq:sigmaPsigma}) sweep the whole Lagrangian Grassmannian.
This proves b).

\smallskip

4. Let  $v\oplus w\in P[\B_1]$, $w\oplus y\in P[\B_2]$, i.e.,
$$
\wh a(w)\, \B_1=\B_1\,\wh a(v),\qquad
\wh a(y)\, \B_2= \B_2\,\wh a(w)
.
$$
Then
$$
\B_2\B_1 \wh a(v)=\B_2\wh a (w) \B_1=\wh a(y) \B_2\B_1
$$
and this proves c).

We omit a proof of d), which is more difficult (see \cite{Ner-book}, Subsect. II.6.5).
\hfill $\square$


\smallskip

{\bf \punct Explicit correspondence. Another description.}
The  proof above  implies
also the following theorem (see \cite{Ner-book}, Theorem II.6.11).

\begin{theorem}
\label{th:spinor-correspondence-2}
Let $P$ satisfy Lemma {\rm\ref{l:potapov-2}}.
Then the corresponding Berezin  operator $\frB[P]$
has the kernel
$$
\exp\Bigl\{\frac12 \begin{pmatrix}\xi&\ov \eta \end{pmatrix}
\begin{pmatrix} A&-B\\ B^t& -C\end{pmatrix}
\begin{pmatrix}\xi^t\\ \ov \eta^t \end{pmatrix}
\Bigr\}
.$$
\end{theorem}


\section{A survey of the oscillator representation.\\
 The category of Gaussian integral
operators}

\label{s:symplectic}

\COUNTERS

In Subsections \ref{ss:fock-space}--\ref{ss:fock-operators}, we
define the bosonic Fock space. In Subsections
\ref{ss:gaussian-operators}--\ref{ss:product-fock}, we introduce
Gaussian integral operators. The algebraic structure of the category
of Gaussian integral operators is described in Subsections
\ref{ss:category-Sp}--\ref{ss:relation-from-operator}.
For a detailed exposition, see \cite{Ner-book}, Chapter  4, or \cite{Ner-lectures}, Sect. 5.1-5.2


\smallskip

{\bf\punct Fock space.%
\label{ss:fock-space}}
 Denote by $d\lambda(z)$ the Lebesgue measure
on $\C^n$ normalized as
$$
d\lambda(z):=\pi^{-n} dx_1\dots dx_n\, dy_1\dots dy_n,\qquad
 \text{where $z_j=x_j+iy_j$.}
$$
 The bosonic Fock space $\bF_n$ is the
space of entire functions on $\C^n$ satisfying the condition
$$
\int_{\C^n} |f(z)|^2 e^{-|z|^2}\, d\lambda(z)<\infty
.
$$
 We define the inner product in $\bF_n$ by the formula
$$
\la f,g\ra:=\int_{\C^n} f(z)\ov{g(z)} e^{-|z|^2}\, d\lambda(z)
.
$$

\begin{theorem}
\label{th:fock-complete}
The space $\bF_n$ is complete, i.e., it is a Hilbert space.
\end{theorem}

\begin{proposition}
\label{pr:fock-basis}
The monomials $z_1^{k_1}\dots z_n^{k_n}$ are pairwise orthogonal and
$$
\|z_1^{k_1}\dots z_n^{k_n}\|^2=\prod k_j!
$$
\end{proposition}


{\bf\punct Operators.%
\label{ss:fock-operators}}

\begin{theorem}
\label{th:operator-symbol}
For each bounded operator $A:\bF_n\to\bF_m$,
there is a function $K(z,\ov u)$
on $\C^m\oplus\C^n$ holomorphic in $z\in\C^m$ and
antiholomorphic
in $u\in \C^n$ such that
$$
Af(z)=\int_{\C^n} K(z,\ov
u)\,f(u)\,e^{-|u|^2}\,d\lambda(u)
$$
{\rm(}the integral  absolutely converges for all $f${\rm)}.
\end{theorem}


{\bf \punct Gaussian operators.%
\label{ss:gaussian-operators}}
Fix $m$, $n=0,1,2,\dots$
Let
$S=\begin{pmatrix}
K&L\\ L^t& M
\end{pmatrix}
$ be a symmetric $(m+n)\times(m+n)$-matrix, i.e.,
$S=S^t$.
A Gaussian operator
$$
\B[S]=\B  \begin{bmatrix}
K&L\\ L^t& M
\end{bmatrix} : \F_n\to \F_m
$$
is defined by
$$
\B \begin{bmatrix}
K&L\\ L^t& M
\end{bmatrix} f(z)=\int_{\C^n}
\exp\Bigl\{\frac 12
\begin{pmatrix}z&\ov u\end{pmatrix}
 \begin{pmatrix} K&L\\ L^t&M \end{pmatrix}
\begin{pmatrix}z^t\\ \ov u^t\end{pmatrix}\Bigr\}
f(u)\, e^{-|u|^2}\,d\lambda(u)
,
$$
where
$$z:=\begin{pmatrix}z_1&\dots&z_n\end{pmatrix},
\qquad
\ov u:=\begin{pmatrix}\ov u_1&\dots&\ov u_n\end{pmatrix}
$$
are row-matrices.

For an operator $A$ from the standard Euclidean space $\C^l$ to the standard Euclidean space $\C^k$ denote by $\|A\|$ the
operator norm of $A$,
$$
\|A\|=\max_{v\in \C^l,\,\|v\|=1}\|Av\|=\Bigl\{\text{the maximal eigenvalue of $A^*A$} \Bigr\}^{1/2}.
$$

\begin{theorem}
\label{th:olshanski}
{\rm(}G. I. Olshanski{\rm)}
An operator $\B[S]$ is bounded
if and only if%
\footnote{The condition 1 implies $\|K\|\le 1$, $\|M\|\le 1$, the condition 2 strengthen these inqualities.}

\smallskip

{\rm 1.} $\|S\|\le 1$,

\smallskip

{\rm 2.} $\|K\|<1$, $\|M\|<1$.
\end{theorem}


{\bf\punct Product formula.%
\label{ss:product-fock}}

\begin{theorem}
\label{th:product-olshanski}
Let
$$\B[S_1]=\B\begin{bmatrix} P&Q\\Q^t& R \end{bmatrix}:
\F_n\to \F_m,\qquad
\qquad
\B[S_2]=\B\begin{bmatrix} K&L\\L^t& M \end{bmatrix}:
\F_m\to \F_k
$$
be bounded Gaussian operators. Then their product
is
$$
\det(1-MP)^{-1/2}\,\, \B[S_2*S_1]
,
$$
 where
$S_2*S_1$ is given by
\begin{multline}
\begin{pmatrix} K&L\\L^t& M \end{pmatrix}
*
\begin{pmatrix} P&Q\\Q^t& R \end{pmatrix}
=\\=
\begin{pmatrix}
K+LP(1-MP)^{-1}L^t & L(1-PM)^{-1}Q\\
Q^t(1-MP)^{-1} L^t& R+Q^t(1-MQ)^{-1}MQ
\end{pmatrix}
\label{eq:strange}
.
\end{multline}
\end{theorem}

\begin{theorem}
\label{th:weil-representation}
Denote by $G_n$ the set of unitary $(n+n)\times(n+n)$
symmetric matrices $\begin{pmatrix}K&L\\L^t&M\end{pmatrix}$
satisfying $\|K\|<1$, $\|M\|<1$. Then $G_n$ is closed
with respect to the $*$-multiplication and is isomorphic
to the group $\Sp(2n,\R)$.
\end{theorem}

Observe that   formula (\ref{eq:strange}) almost
coincides with  formula (\ref{eq:strange-berezin}).
Again, formula (\ref{eq:strange}) hides a product of linear relations.

\smallskip

First, we define analogues of the spaces
$\cV_{2n}$ from the previous section.

\smallskip


{\bf\punct Complexification
 of a linear space with bilinear form.%
\label{ss:objects-Sp}}
Denote by $\cW_{2n}$ the space $\C^{2n}=\C^n\oplus \C^n$.
Let us denote its elements by
$v=\begin{pmatrix} v^+ & v^-\end{pmatrix}$.
We equip $\cW_{2n}$ with two forms

\smallskip

--- the skew-symmetric bilinear form
\begin{equation}
\Lambda(v,w)=\Lambda_n(v,w):=\begin{pmatrix} v^+ & v^-\end{pmatrix}
\begin{pmatrix}
0&1_n\\ -1_n&0
\end{pmatrix}
\begin{pmatrix} (w^+)^t \\ (w^-)^t\end{pmatrix}
=v^+ (w^-)^t-v^- (w^+)^t
;
\label{eq:form-Lambda}
\end{equation}

--- the indefinite Hermitian form
\begin{equation}
M(v,w)=M_n(v,w):=\begin{pmatrix} v^+ & v^-\end{pmatrix}
\begin{pmatrix}
1_n&0\\ 0&-1_n
\end{pmatrix}
\begin{pmatrix} (\ov v^+)^t \\ (\ov w^{\,-})^t\end{pmatrix}
=v^+(\ov v^+)^t- w^-(\ov w^{\,-})^t
.
\label{eq:form-M}
\end{equation}

{\sc Remark.}
 Let us  explain the origin of the definition.
Consider a real space $\R^{2n}$ equipped with
 a nondegenerate
skew-symmetric bilinear form $\{\cdot,\cdot\}$.
 Consider the space $\C^{2n}\supset\R^{2n}$.
We can extend $\{\cdot,\cdot\}$ to $\C^{2n}$ in
the following two ways.
First, we can  extend it as a bilinear form,
$$
\wt\Lambda(x+iy,x'+iy'):=\{ x,x'\}-
 \{y ,y' \}+ i\bigl(\{ x,y' \}+ \{x' ,y \}\bigr)
.
$$
Next, we  extend $\{\cdot,\cdot\}$ as a sesquilinear form
$$
\wt M(x+iy,x'+iy') := \{x ,x' \}+ \{y ,y' \} +i\bigl( \{y ,x' \}- \{x ,y' \} \bigr)
.
$$
This form is anti-Hermitian, $\wt M(u,z)=-\ov{\wt M(z,u)}$, it is more convenient to pass to 
a bilinear form
$\Lambda:=i\wt\Lambda$, and a Hermitian form $M:=i\wt M$. 
Thus, we get a space endowed with two forms, skew symmetric and Hermitian.
Denoting by $e_l$ the standard basis in $\C^{2n}$ and passing to a new basis
$(e_j+e_{n+j})/\sqrt 2$,  $(e_j-e_{n+j})/\sqrt 2$, where $j\le n$, we arrive at expressions 
\eqref{eq:form-Lambda}, \eqref{eq:form-M}.
\hfill $\boxtimes$


\smallskip

{\bf\punct  The category $\bSp$.%
\label{ss:category-Sp}} 
Objects of the category $\bSp$ are the spaces
$\cW_{2n}$, where $n=0$, $1$, $2$, \dots.

 For $m$, $n=0$, 1, 2, \dots, we equip the direct sum
$\cW_{2n}\oplus\cW_{2m}$ with two forms,
\begin{align}
\Lambda^\ominus(v\oplus w,v'\oplus w'):=
\Lambda_n(v,v')-\Lambda_m(w,w'),
\nonumber
\\
M^\ominus(v\oplus w,v'\oplus w'):=
M_n(v,v')-M_m(w,w').
\label{eq:M-ominus}
\end{align}

A morphism $ \cW_{2n}\to\cW_{2m}$
is a linear relation $P:\cW_{2n}\arr\cW_{2m}$
satisfying the following conditions.

\smallskip

1. $P$ is $\Lambda^\ominus$-Lagrangian.

\smallskip

2. The form $M^\ominus$ is non-positive on $P$.

\smallskip

3. The form $M_{\cW_{2n}}$ is strictly negative on $\ker P$
and the form  $M_{\cW_{2m}}$
is strictly positive on $\indef P$.

\smallskip

A product of morphisms is the usual product of linear relations%
\footnote{Under the condition 3, for morphisms $P:\cW_{2n}\arr\cW_{2m}$, $Q:\cW_{2m}\arr\cW_{2k}$, we have $\im P\cap \ker Q=0$. For this reason
 $\nul$ does not appear in the category $\bSp$.}.

 \begin{observation}
\label{obs:automorphisms-sp}
The group of automorphisms of $\cW_{2n}$ is the real
symplectic group $\Sp(2n,\R)$.
\end{observation}

This follows from the remark given in the previous subsection.
The group of operators preserving  both the forms $\wt\Lambda$,
$\wt M$ on $\C^{2n}$ preserves also the real subspace $\R^{2n}$.


\smallskip

{\bf\punct Construction of Gaussian operators from  linear
relations.%
\label{ss:operator-from-relation}}
 Recall that $\cW_{2n}$ is $\C^{2n}=\C^n\oplus\C^n$.
Denote this decomposition by
$$
\cW_{2n}=\cW_n^+\oplus\cW_n^-
.
$$
Represent a linear relation $P$ as the graph of an operator
$$
S=S(P):\cW_m^-\oplus \cW_n^+  \to
\cW_m^+\oplus \cW_n^-
.
$$
This is possible, because $M^\ominus$ is negative semidefinite on the subspace  $P$ and  is strictly positive on the subspace
$\cW_m^+\oplus \cW_n^-$;
therefore $P\cap(\cW_m^+\oplus \cW_n^-)=0$.

\begin{proposition}
\label{pr:potapov-transform}
A matrix $S$ has the form $S(P)$ if and only if it is symmetric and satisfies
the Olshanski conditions from Theorem {\rm\ref{th:olshanski}}.
\end{proposition}

\begin{theorem}
\label{th:weil-functor}
For each morphisms
$$
P:\cW_{2n}\arr\cW_{2m},\, \qquad Q:\cW_{2m}\arr\cW_{2k}
,
$$
the corresponding Gaussian operators
$$
\B\bigl[S(P)\bigr]:\F_n\to\F_m,\qquad \B\bigl[S(Q)\bigr]:\F_m\to\F_k
$$
satisfy
$$
\B\bigl[S(Q)\bigr]\, \B\bigl[S(P)\bigr]=\lambda(Q,P)\, \B\bigl[S(QP)\bigr]
,
$$
where $QP$ is the product of linear relations and
$\lambda(Q,P)$ is a nonzero scalar.
\end{theorem}

As  formulated, the theorem can be proved
by  direct force.


\smallskip

{\bf\punct Construction of linear relations
from Gaussian operators.%
\label{ss:relation-from-operator}}
For
$$
\begin{pmatrix}v_1^+&\dots &v_n^+& v_1^-&\dots& v_n^-  \end{pmatrix}
\in \cV_{2n}
,
$$
we define the differential operator
 (a {\it creation-annihilation operator})
$$
\wh a(v) f(z)=
\Bigl(\sum_j v_j^+ z_j +
\sum_j v_j^- \frac{\partial}{\partial z_j}  \Bigr)
\,f(z)
.
$$

For a given bounded Gaussian operator $\B[S]:\F_n\to\F_m$ we consider the set $P$ of all
$v\oplus w\in\cW_{2n}\oplus \cW_{2m} $
such that
$$
\wh a(w)\,\B[S]=\B[S]\,\wh a(v).
$$

\begin{theorem}
The linear relation $P$
 is a morphism of the category $\bSp$.
\end{theorem}

{\bf\punct Details. An analogue of the Schwartz space.%
\label{ss:schwartz-space}}
We define the  Schwartz--Fock space $\cS\F_n$ as  the subspace in $\F_n$ consisting
of all
$$
f(z)=\sum c_{j_1,\dots, j_n} z^{j_1}\dots z^{j_n}
$$
such that for each $N$
$$
\sup_j |c_{j_1,\dots, j_n}| \prod_k j_k!\, {j_k}^N
<\infty
.$$

\begin{theorem}
\label{th:gauss-schwartz}
The subspace $\cS\F_n$ is a common invariant domain for all
 Gaussian
bounded operators
and for all  creation-annihilation operators.
\end {theorem}

See \cite{Ner-lectures}, Subsect. 4.2.4 and Theorem 5.1.5.


\smallskip

{\bf\punct Details. The Olshanski
 semigroup $\Gamma\Sp(2n,\R)$.%
\label{ss:olshanski}}
The {\it Olshanski semigroup} $\Gamma\Sp(2n,\R)$
is defined as the subsemigroup in $\Sp(2n,\C)$
consisting of complex matrices $g$ satisfying the condition
$$
g\begin{pmatrix}-1&0\\0&1\end{pmatrix}g^*
-\begin{pmatrix}-1&0\\0&1\end{pmatrix}\le 0
,
$$
where $g^*=\ov g^t$ denotes the adjoint matrix (see \cite{Olsh-semi}).

Equivalently,
 $g\in\Gamma\Sp(2n,\R)$
if and only if
$$
M(ug,ug)\le M(u,u)\q
\text{for all $u\in\C^{2n}$}
.$$

The Olshanski semigroup is a subsemigroup
in the semigroup of endomorphisms of the
object $\cW_{2n}$. For details, see \cite{Ner-lectures}, Sect. 2.7, 3.5.


\section{Gauss--Berezin integrals}

\label{s:integral}

\COUNTERS

Here we discuss  super-analogues of Gaussian integrals.
Actually, the final formulas  are not used, but their structure
is important for us.

Apparently, these integrals are evaluated somewhere,
but I do not know a reference.
 Calculations in the fermionic case
are contained in \cite{SMD}.

\smallskip


{\bf\punct Phantom algebra.%
\label{ss:phantom}}
 {\it Phantom generators} $\fra_1$, $\fra_2$, \dots
are anticommuting variables,
$$
\fra_k \fra_l=-\fra_l\fra_k, \qquad \fra_j^2=0
.
$$
We define a {\it phantom algebra} $\cA$
as the algebra of polynomials in the  variables
$\fra_j$. {\it For the sake of simplicity, we assume
that the number of  variables is infinite}.
We also call elements of $\cA$  {\it phantom constants}.

The phantom algebra has a natural $\Z$-grading by
 degree of  monomials,
$$
\cA=\oplus_{j=0}^\infty \cA_j
.
$$
Therefore, $\cA$ admits a $\Z_2$-grading, namely
$$
\cA_\ev:=\oplus \cA_{2j},
\qquad
\cA_\od:=\oplus \cA_{2j+1}
.
$$
We define the automorphism $\mu\mapsto \mu^{\sigma}$
of $\cA$
by the rule
\begin{equation}
\mu^\sigma=\begin{cases}
\mu \quad &\text{if $\mu$ is even,}
\\
-\mu \quad &\text{if $\mu$ is odd}
\end{cases}
\label{eq:involution}
\end{equation}
(equivalently, $\fra_j^\sigma=-\fra_j$).

The algebra $\cA$ is {\it supercommutative}
in the following sense:
\begin{align*}
\mu\in\cA, \,\nu\in\cA_\ev\q
 \Longrightarrow \q \mu\nu=\nu\mu,
\\
\mu\in\cA_\od, \,\nu\in\cA_\od
\q\Longrightarrow\q \mu\nu=-\nu\mu
.
\end{align*}
Also,
\begin{equation}
\mu\in\cA,\, \nu\in\cA_{\odd}
 \Longrightarrow \q \nu\mu=\mu^\sigma\nu.
\label{eq:sigma!}
\end{equation}

Next, represent $\mu\in\cA$
as $\mu=\sum_{j\ge 0} \mu_j$,
 where $\mu_j\in\cA_j$.
 We define the map
$$
\pia:\cA\to\C
$$
by
 $$
 \pia(\mu)=
\pia\Bigl(\sum\nolimits_{j\ge 0} \mu_j\Bigr)
:=
 \mu_0\in\C
. $$
Evidently,
$$
\pia(\mu_1\mu_2)=\pia(\mu_1)\,\pia(\mu_2), 
\qquad \pia(\mu_1+\mu_2)=\pia(\mu_1)+\pia(\mu_2)
.$$

Take $\phi\in\cA$ such that $\pia(\phi)=0$. Then
$\phi^N=0$ for sufficiently large $N$. Therefore,
\begin{equation}
(1+\phi)^{-1}:=\sum_{n\ge 0} (-\phi)^n
,
\label{eq:inverse}
\end{equation}
actually, the sum is finite. In particular,
if $\pia(\mu)\ne 0$, then $\mu$ is invertible.


\smallskip

{\bf\punct A technical comment.%
\label{ss:super}}
The aim of this paper is a specific construction
and we  use a minimal vocabulary necessary for our aims.
We regard  supergroups $\GL(p|q)$ and $\OSp(2p|2q)$ as groups of matrices  over the algebra $\cA$, representations of supergroups 
 are defined in  modules  over $\cA$.

The more common point of view%
\footnote{On comparison and criticism of different definitions of super-objects, see, e.g., \cite{Mol}, \cite{Sch}.
Our approach follows DeWitt's book \cite{DeW},  we use only linear algebra and integral oprators (without analysis on manifolds). 
Wider generality discussed in this subsection is similar to Berezin--Kats \cite{BK}.}
is to consider supergroups as functors from the category of
supercommutative algebras to the category of groups. Our constructions require a restriction
of a class of supercommutative algebras%
\footnote{We need exponentials of even elements \eqref{eq:exp},
$\C$-valued integrals \eqref{eq:integration}, and
inverses \eqref{eq:inverse}, \eqref{eq:GB-integral}. We also must justify
the calculation in Theorem \ref{th:GB-integral} and the proof of Lemma \ref{l:intersection}.}.
Precisely, the algebra $\cA$  can be  replaced%
\footnote{with some modifications in proofs.} by an arbitrary
finitely or countably generated supercommutative
algebra $\cA'$ over $\C$ satisfying the following properties:

\sm

--- there is a non-zero homomorphism $\pi_\downarrow:\cA'\to \C$,
denote $I:=\ker \pi_\downarrow$;

\sm

--- any finitely generated subalgebra in $I$ is nilpotent.

\sm

A pass to algebras equipped with the Krull topology make situation more flexible.
We can consider supercommutative algebras $\cA''$ satisfying the following properties: 

---  there is a non-zero homomorphism $\pi_\downarrow:\cA''\to \C$,
denote $I:=\ker \pi_\downarrow$;

\sm

--- $I$ is finitely or countably generated, and $\cap_n I^n=0$;

\sm

--- $\cA''$ is complete with respect to $I$-adic topology (Krull topology, see
 \cite{Eis}, Subsect. 7.7), i.e., for any sequence $a_n\in \cA''/I^n$ such that natural maps
 $\cA''/I^n\to \cA''/I^{n-1}$ send $a_n$ to $a_{n-1}$, there is $a\in \cA''$
 whose image in each $\cA''/I^n$ is $a_n$.

\sm

All our constructions  below are functorial with respect to  algebras of such types.

\sm

Below we use terms '{\it supergroups}' and {\it super-Grassmannians} for objects defined over the
algebra $\cA$ keeping in mind that this can be translated into the functorial language.

\smallskip


{\bf\punct Berezinian.%
\label{ss:berezinian}}
 Let $\begin{pmatrix}P&Q\\R&T\end{pmatrix}$
be a block $(p+q)\times (p+q)$-matrix,
let $P$, $T$ be composed of
even phantom constants, and $Q$, $R$ be composed
of odd phantom constants. Then the Berezinian
(or Berezin determinant)
is%
\footnote{The usual determinant
of a block complex matrix $\begin{pmatrix}P&Q\\R&T\end{pmatrix}$ is
$\det P \cdot\det(T-QP^{-1}R)$.}
$$
\mathrm{ber} \begin{pmatrix}P&Q\\R&T\end{pmatrix}
:=(\det P)^{-1} \cdot \det(T-QP^{-1}R)
.
$$
We note  that $P$ and $T-QP^{-1}R$ are composed
of elements of the commutative algebra
$\cA_{\even}$, therefore their determinants are well-defined.
The Berezinian satisfies the multiplicative
 property of the usual determinant
$$
\mathrm{ber}( A)\,\mathrm{ber} (B)
=\mathrm{ber} (AB)
.
$$


\smallskip

{\bf\punct Functions.%
\label{ss:function}} We consider 3 types of variables:

\smallskip

--- real or complex (bosonic) variables, we denote them by $x_i$,
$y_j$
(if they are real) and $z_i$, $u_j$ (if they are complex);

\smallskip

--- Grassmann (fermionic) variables, we denote them by $\xi_i$, $\eta_j$
or $\ov\eta_j$;

\smallskip

--- phantom generators $\fra_j$ as above.

\smallskip

Bosonic variables $x_l$
commute with the fermionic variables $\xi_j$ and phantom
constants $\mu\in\cA$. We also assume that
the fermionic variables $\xi_j$
and the phantom generators $\fra_l$ anticommute,
$$
\xi_j \fra_l=-\fra_l\xi_j
.
$$

Fix a collection of bosonic variables
$x_1$, \dots, $x_p$ and a collection
of fermionic variables $\xi_1$,\dots, $\xi_q$.
A {\it function} is a
sum of the form
\begin{multline}
f(x,\xi):=
\sum_{m\ge 0}\,\,
\sum_{0<i_1<\dots<i_k\le q;\, j_1<\dots<j_m} 
 h_{i_1,\dots,i_k; j_1,\dots, j_m}(x_1,\dots, x_p)
 \times\\\times
 \fra_{j_1}\dots \fra_{j_m} \,
\xi_{i_1}\dots\xi_{i_k},
\label{eq:function}
\end{multline}
where
 $h$ are smooth functions
 of $x\in\R^p$. We also  write such expressions in the form
$$
f(x,\xi)= \sum_{I,J} h_{I,J}(x) \fra^J \xi^I
$$ 
 keeping in the mind that $I$ ranges in collections $0<i_1<\dots<i_k\le q$
 and $J$ in collections $j_1<\dots<j_m$.

\smallskip

We say, that a function $f$ is {\it even}
(respectively {\it odd})
if it is an even expression in
the total collection $\xi_i$, $\fra_k$.
By $f^\sigma(x,\xi)$ we denote the function obtained from 
$f(x,\xi)$ by the substitution $\fra_j\mapsto -\fra_j$ for all
$j$.

\smallskip

{\sc Remark.} Formally,
the fermionic variables and the phantom constants have equal rights
in our definition.
However, below their roles are different: $\xi_j$
serve as variables and elements of $\cA$ serve as constants
(see (\ref{gauss-berezin-expression})).
\hfill $\boxtimes$

\sm

{\sc Remark.} There are three ways  to understand expressions
\eqref{eq:function}. We consider them as finite sums, but it is possible 
to consider them as arbitrary formal series in variables $\fra_j$.
We also can consider formal series such that each $m$-th summand depends only on finitely many of $\fra_j$ (considerations below
survive in these cases after minor modifications).
\hfill $\boxtimes$

\smallskip

For a given $f$, we define the function
$$
\pia(f):= \sum_I h_{I,\emptyset}(x)\, \xi^I\,\in\, C^\infty(\R^p)\otimes \Lambda_q
.
$$


{\bf\punct Integral.%
\label{ss:super-integral}} Now, we define the  symbols
$$
\int f(x,\xi)\,dx\, \qquad \int f(x,\xi)\,d\xi,
\qquad
\int f(x,\xi)\,dx\,d\xi
.
$$
The integration with respect to $x$ is the usual
 termwise integration in (\ref{eq:function}),
\begin{equation}
\int_{\R^p} f(x,\xi)\,dx:=
\sum_{I,J}\Bigl(\int h_{I,J}(x)\,dx\Bigr) \fra^J \xi^I.
\label{eq:integration}
\end{equation}

The integration with respect to $\xi$ is the usual termwise Berezin
integral,
$$
\int f(x,\xi)\,d\xi
:=
\sum_{J}  h_{L J}(x)
\fra^J , \qquad \text{where $L=\{1,2,\dots,q\}$}. 
$$


{\bf\punct  Exponential.%
\label{super-exponent}}
Let $f(x,\xi)$ be an {\it even} expression in $\xi$, $\fra$,
i.e., $f(x,\xi)=f(x,-\xi)^{\sigma}$.
 We define its exponential
as usual, it satisfies the usual properties.
Namely,
\begin{equation}
\exp\{f(x,\xi)\}:=\sum_{n=0}^\infty \,\frac1{n!} f(x,\xi)^n
\label{eq:exp}
.
\end{equation}

Since $f_1$, $f_2$ are even, we have $f_1 f_2=f_2f_1$. Therefore, the  identity
$$
\exp\{f_1+f_2\}=\exp\{f_1\} \exp\{f_2\}
$$
holds.

\begin{observation}
The series {\rm(\ref{eq:exp})} converges.
\end{observation}

Indeed,
\begin{multline*}
\exp\bigl\{   f(x,\xi)\bigr\}
=
\exp\bigl\{
h_{\emptyset\emptyset}(x)
 \bigr\}
\prod_{(I,J)\ne (\emptyset,\emptyset)}
 \exp\bigl\{ h_{I,J}(x) \fra^I \xi^J \bigr\}
=\\=
\exp\bigl\{ h_{\emptyset\emptyset}(x) \bigr\}
\prod_{(I,J)\ne(\emptyset,\emptyset)}
(1+ h_{I,J}(x) \fra^I \xi^J)
.
\end{multline*}
Opening brackets, we get
a polynomial in $\fra_j$, $\xi_k$.


\smallskip

{\bf\punct Gauss--Berezin integrals. A special case.%
\label{ss:GB-integral}}
Take $p$ real variables $x_i$ and $q$ Grassmann
variables $ \xi_j $.
Consider the expression
\begin{multline}
I=
\iint \exp\Bigl\{\frac 12
\begin{pmatrix} x&\xi\end{pmatrix}
\begin{pmatrix}
A&B\\ -B^t& C
\end{pmatrix}
\begin{pmatrix} x^t\\ \xi^t\end{pmatrix}
\Bigr\}\,dx\,d\xi
=\\=
\iint
\exp\Bigl\{
\frac 12 \sum_{ij} a_{ij} x_i x_j
+\sum_{ik}b_{ik}x_i\xi_k
+\frac 12 \sum_{kl} c_{kl} \xi_k\xi_l
\Bigr\}
\,dx\,d\xi
\label{gauss-berezin-expression}
.
\end{multline}

The notation $^t$ denotes the transpose as above, also
$\begin{pmatrix} x&\xi\end{pmatrix}$
 is the row-matrix
$$
\begin{pmatrix} x&\xi\end{pmatrix}=
\begin{pmatrix} x_1&\dots & x_p& \xi_1 & \dots &\xi_q\end{pmatrix}
.
$$
The matrices $A$ and $C$
 are composed of even phantom constants,
$A$ is symmetric, $C$ is skew-symmetric,
and $B$ is a matrix composed of odd phantom constants.%
\footnote{The argument of the exponential must be even in $\xi$, $\fra$.
This imposes the constraints
 of parity for $A$, $B$, $C$.
The symmetry conditions for $A$, $B$, $C$ are the natural
conditions for coefficients of a quadratic form in $x$, $\xi$.%
\label{fo:}}

\begin{observation}
The integral  converges if and only if
the matrix $\Re \pia(A)$ is negative definite.
\end{observation}

 Indeed,
the integrand $\exp\{\dots\}$
is a finite sum
of the form
$$
\exp\Bigl\{\frac 12
\sum_{ij}\pia(a_{ij})x_i x_j\Bigr\}
\sum_{i_1<\dots<i_k}\,\,\sum_{j_1<\dots<j_l}
P_{i_1,\dots,i_k;\,j_1,\dots,j_l}(x)
\xi_{i_1}\dots\xi_{i_k}
\fra_{j_1}\dots\fra_{j_l}
,
$$
where $P_{\dots}(x)$ are polynomials.
Under the condition  $\Re \pia(A)<0$,
 a term-wise integration is possible.
 \hfill $\square$


\smallskip

{\bf\punct Evaluation of the Gauss--Berezin integral.%
\label{ss:evaluation-GB}} Let us evaluate integral \eqref{gauss-berezin-expression}. 

\begin{theorem}
\label{th:GB-integral}
Let $\Re \pi_{\downarrow}(A)<0$. Then
\begin{equation}
I=
\begin{cases}
(2\pi)^{p/2}\,\det(-A)^{-1/2}\,\Pfaff(C+B^tA^{-1}B)
\quad&\text{if $q$ is even}
,
\\
0,\quad&\text{otherwise}
.
\end{cases}
\label{eq:GB-integral}
\end{equation}
\end{theorem}

Recall that $q$ is the number of Grassmann variables.

\smallskip

{\sc Remark.} The matrix $C+B^tA^{-1}B$ is
skew-symmetric and composed of
even phantom constants.
Therefore, the Pfaffian is well defined.
\hfill $\boxtimes$

\smallskip

{\sc Remark.}
Thus, for $q$  even
our expression is a 'hybrid'  of a Pfaffian and a Berezinian,
$$
I^{-2}=-(2\pi)^p \,\mathrm{ber}
\begin{pmatrix}
A&B\\-B^t&C
\end{pmatrix}
.
$$
Similar (but not precisely same) 'hybrid' appeared in \cite{Sergeev} and \cite{SupPf}.
\hfill $\boxtimes$

\sm

{\sc Proof of Theorem \ref{th:GB-integral}.} First, we integrate with respect to $x$,
\begin{multline*}
\exp\Bigl\{\frac 12 \xi C\xi^t\Bigr\}
\int_{\R^p}
\exp\Bigl\{
\frac 12 xAx^t
+xB\xi^t
\Bigr\}
\,dx
=\\=
\exp\Bigl\{\frac 12 \xi C\xi^t\Bigr\}
\times \\ \times
\int_{\R^p}
\exp\Bigl\{
\frac 12 (x-\xi B^t A^{-1})
A
(x^t+A^{-1}B\xi^t )\Bigr\}
\exp\Bigl\{
\frac 12\xi B^t A^{-1} B\xi^t
\Bigr\}
\,dx
.\end{multline*}

We substitute
\begin{equation}
y:=x-\xi B^tA^{-1}
\label{eq:substitution}
,\end{equation}
  get
$$
\exp\Bigl\{\frac 12 \xi C\xi^t\Bigr\}\cdot
\exp
\Bigl\{
\frac 12 \xi B^t A^{-1} B\xi^t\Bigr\}
 \int \exp\Bigl\{\frac 12 y A y^t\Bigr\}\,dy
 ,
$$
and arrive at the usual Gaussian integral
(\ref{eq:gaint}).

Integrating
the result, we get
$$
\det(-A)^{-1/2}(2\pi)^{p/2}
\int
\exp\Bigl\{\frac 12 \xi (C+ B^t A^{-1} B)\xi^t\Bigr\}
\,d\xi
,
$$
and arrive at the Pfaffian.

We must justify the substitution
(\ref{eq:substitution}).
Let $\Phi$ be a function on $\R^p$ of Schwartz class, let
$\nu$ be an even expression in $\fra$, $\xi$, assume
that the constant term of $\nu$ is 0.
 Then
$$
\int_{\R^p} \Phi(x+\nu)\,dx=\int_{\R^p} \Phi(x)\,dx
.
$$
Indeed,
$$
\Phi(x+\nu):=\sum_{j=0}^\infty \frac 1{j!}
\nu^j\frac {d^j}{dx^j} \Phi(x).
$$
Actually, the summation is finite.
A termwise integration
with respect to $x$ gives zero
for all $j\ne 0$.
\hfill $\square$


\smallskip

{\bf\punct  Grassmann Gaussian integral.%
\label{ss:grassmann-gauss}}

\begin{observation}
\label{obs:grassmann-gauss-integral}
Let  $D$ be a  {\sf complex} skew-symmetric
matrix of size $N$, let $\xi_k$, $\zeta_k$
be  Grassmann variables.
Then the integral
$$
\int \exp\Bigl\{ \frac 12 \xi D \xi^t+ \xi \zeta^t\Bigr\}
\,
d\xi
$$
can be represented in the form
$$
s\cdot \prod_{j=1}^m \bigl( \zeta h_j^t \bigr)
\cdot
\exp\Bigl\{\frac 12 \zeta Q\zeta^t\Bigr\}
,
$$
where $s\in\C$,  $Q$ is a skew-symmetric matrix,
and
$h_j$ are row-matrices,
$m\le N$, $N-m$ is even.
\end{observation}

Indeed, one can
find a linear substitution
$
\xi=\eta S
$
such that%
\footnote{In other words, one can reduce a skew-symmetric
matrix over $\C$ to a canonical form.}
$$
\xi D\xi^t=\sum_{j=1}^\gamma \eta_{2j-1}\eta_{2j1}
.$$
Then the integral can be reduced
to
$$
\det S\int \exp\Bigl \{\sum_{j=1}^\gamma \eta_{2j-1}\eta_{2j}
+
\sum_{k=1}^{N} \eta_k\nu_k\Bigr\}\,\, d\eta
,
$$
where $\nu_j$ are certain linear expressions in $\zeta_l$.
So we get
\begin{multline*}
\det S\cdot \prod_{j=1}^\gamma
 \int \exp\bigl \{ \eta_{2j-1}\eta_{2j}+\eta_{2j-1}\nu_{2j-1}+
\eta_{2j}\nu_{2j}\bigl\}\,d\eta_{2j-1}\, \,d\eta_{2j}
\,\times
\\
\times
\int\exp\Bigl\{ \sum_{j=2\gamma+1}^N \eta_j\nu_j\Bigr\}\,\, d\eta_{2\gamma+1}\dots d\eta_N
=\\=
\pm \det S\cdot \exp\Bigl\{-\sum_{j=1}^\gamma \nu_ {2j}\nu_{2j+1}\Bigr\}
\prod_{k= 2\gamma+1}^N \nu_k
.
\end{multline*}
Recall that $\nu_j$ are certain
 linear expressions%
 \footnote{A product of functions $\nu_m$ is canonically defined up to a constant
factor, equivalently a linear span of functions $\nu_m$ is canonically defined
(this sentence is a rephrasing of the Pl\"ucker embedding of a Grassmannian into
an exterior algebra).}
  in $\zeta_l$.
\hfill $\square$.


\smallskip

{\bf\punct
 More general Gauss--Berezin integrals.%
\label{ss:general-GB}}
Consider an expression
\begin{multline}
J=
\iint \exp\Bigl\{\frac 12
\begin{pmatrix} x&\xi\end{pmatrix}
\begin{pmatrix}
A&B\\ -B^t& C
\end{pmatrix}
\begin{pmatrix} x^t\\ \xi^t\end{pmatrix}
+x h^t+ \xi g^t
\Bigr\}\,dx\,d\xi
=\\=
\iint
\exp\Bigl\{
\frac 12 \sum_{ij} a_{ij} x_i x_j
+\sum_{ik}b_{ik}x_i\xi_k
+\frac 12 \sum_{kl} c_{kl} \xi_k\xi_l
+\\+\sum_j h_j x_j -\sum_k g_k \xi_k
\Bigr\}
\,dx\,d\xi
\label{gauss-berezin-expression-2}
,\end{multline}
here $A$, $B$, $C$ are  as above and $h^t$, $g^t$
are column-vectors, $h_j\in\cA_\ev$, $g_k\in\cA_\od$.

\smallskip

We propose two ways to  evaluate of this integral.

\smallskip


{\bf \punct The first way to evaluate.%
\label{ss:first-way}}
 Substituting
$$
\begin{pmatrix} y & \eta\end{pmatrix}
=\begin{pmatrix}
x & \xi
\end{pmatrix}
+
\begin{pmatrix}
h & g
\end{pmatrix}
\begin{pmatrix}
A&B\\ -B^t& C
\end{pmatrix}^{-1}
,
$$
we get
\begin{multline*}
\exp\Bigl\{ \frac 12
\begin{pmatrix}
h& g
\end{pmatrix}
\begin{pmatrix} A&B\\-B^t&C\end{pmatrix}^{-1}
\begin{pmatrix}
h^t\\ g^t
\end{pmatrix}
 \Bigr\}
\times \\ \times
\iint \exp\Bigl\{\frac 12
\begin{pmatrix} y&\eta\end{pmatrix}
\begin{pmatrix}
A&B\\ -B^t& C
\end{pmatrix}
\begin{pmatrix} y^t\\ \eta^t\end{pmatrix}
\Bigr\}\,dy\,d\eta
\end{multline*}
and arrive at Gauss--Berezin integral
 (\ref{gauss-berezin-expression}) evaluated above.

This way is not perfect, because it uses an inversion
of a matrix $\begin{pmatrix} A&B\\-B^t&C\end{pmatrix}$.

\begin{observation}
\label{obs:invertibility-2}
A matrix
$\begin{pmatrix} A&B\\-B^t&C\end{pmatrix}$
 is invertible if and only if $A$ and $C$ are invertible.
\end{observation}

The necessity  is  evident;
to prove the sufficiency, we note that the matrix
$$T:=\begin{pmatrix} A^{-1}&0\\0& C^{-1}\end{pmatrix}
\begin{pmatrix}
A&B\\ -B^t& C
\end{pmatrix}-
\begin{pmatrix}1&0\\0&1\end{pmatrix}
$$
is composed of nilpotent elements of $\cA$.
 We
write out $(1+T)^{-1}=1-T+T^2-\dots$,
and therefore our initial matrix
$\begin{pmatrix} A&0\\0& C\end{pmatrix} (1+T)$ is invertible.
\hfill $\square$

\smallskip

The matrix $A$ is invertible, because $\Re A<0$.

\smallskip

But the matrix
$C$ is  skew-symmetric.

\smallskip

--- If $q$ is even, then
a $q\times q$ skew-symmetric matrix $C$ in  general
 position is invertible. For noninvertible $C$,
 we have a chance to remove uncertainty. This way leads to
 an expression of the form (\ref{eq:general})
 obtained below.

\smallskip

---  If  $q$ is odd, then
$C$ is  non-invertible; our approach is not suitable.

\smallskip


{\bf\punct The second way to evaluate
 of Gauss--Berezin integrals.%
\label{ss:second-way}} First, we integrate with respect to $x$,
\begin{multline*}
\exp\Bigl\{ \frac 12 \xi C\xi^t + \xi g^t \Bigr\}
\int_{\R^p} \exp\Bigl\{\frac 12 xAx^t+ xB\xi^t+ x h^t\Bigr\}\,dx
=\\=
\exp\Bigl\{ \frac 12 \xi C\xi^t + \xi g^t \Bigr\}
\exp\Bigl\{\frac 12 (h-\xi B^tA^{-1}) A (h^t+ A^{-1}B\xi^t)\Bigr\}
\times\\\times
\int_{\R^p} \exp\Bigl\{\frac 12 (x+h-\xi B^t)
 A (x^t+h^t+B\xi^t)\Bigr\}\,dx
.
\end{multline*}
Substituting $y=x+h-\xi B^t A^{-1}$
and integrating with respect to $y$, we get
$$
(2\pi)^{p/2} \det(- A)^{-1/2}
\exp\Bigl\{\frac 12 (h-\xi B^tA^{-1}) A (h^t+ A^{-1}B\xi^t)+
\frac12 \xi C \xi^t +\xi g^t\Bigr\}
.
$$

Next, we must integrate with respect to $\xi$,
our integral has a  form
\begin{equation}
\int \exp\Bigl\{ \frac 12 \xi D \xi^t+ \xi r^t\Bigr\}
\,
d\xi
\label{eq:xiDxi}
,
\end{equation}
where a matrix $D$ is composed of even
phantom constants and a vector $r$ is odd.

\smallskip

{\it If $D$ is invertible,} we shift the argument again
$\eta^t:=\xi^t+D^{-1} r^t$ and get
$$
\exp\Bigl\{\frac 12 r D^{-1} r^t\Bigr\}
\int \exp\Bigl\{\frac 12 \eta D \eta^t\Bigr\}\,d\eta
,
$$
the last integral  is a Pfaffian. This way is equivalent to the approach
discussed in the previous subsection.

\smallskip

{\it Now, consider an arbitrary  $D$.}
The calculation
of Subsection \ref{ss:grassmann-gauss}
does not survive%
\footnote{Let $D$ be a skew-symmetric matrix over
 $\cA_{\even}$.
If $\pia(D)$ is degenerate, then we can not reduce $D$
 to a normal form.}.

However, we can write (\ref{eq:xiDxi})
explicitly as follows. For any subset 
$$I:\,i_1<\dots<i_{2k}$$
in $\{1,\dots,q\}$ we
consider
the complementary subset 
$$J:j_1<\dots<j_{q-2k}.$$
 Define the constant
 $\sigma(I)=\pm 1$ as follows
$$
\bigl( \xi_{i_1}\xi_{i_2}\dots\xi_{i_{2k}} \bigr)
\bigl( \xi_{j_1}\xi_{j_2}\dots\xi_{j_{q-2k}}
 \bigr)=\sigma(I)\,\xi_1\xi_2\dots\xi_q
.
$$

Evidently,
\begin{multline}
\int \exp\Bigl\{ \frac 12 \xi D \xi^t+ \xi r^t\Bigr\}
\,
d\xi
=\\=
\sum_I \sigma(I) \Pfaff\begin{pmatrix}
0&d_{i_1 i_2}&\dots & d_{i_1 i_{2k}}
\\
d_{i_2 i_1}&0&\dots & d_{i_2 i_{2k}}
\\
\vdots & \vdots & \ddots & \vdots
\\
d_{i_{2k} i_1}& d_{i_2 i_{2k}}&\dots & 0
\end{pmatrix}
r_{j_1}\dots r_{j_{q-2k}}
\label{eq:general}
.
\end{multline}
Recall that $d_{pq}\in \cA_\ev$
 and $r\in\cA_\od$.


\section{Gauss--Berezin integral operators}

\label{s:operators}

\COUNTERS

Here we define super hybrids of Gaussian
operators and Berezin operators.


\smallskip

{\bf\punct Fock--Berezin spaces.%
\label{ss:fock-berezin}}
Fix $p$, $q=0$, $1$, $2$, \dots.
Let $z_1$, \dots, $z_p$ be complex variables,
 $\xi_1$, \dots, $\xi_q$
 be Grassmann variables.
We consider expressions
$$
f(z,\xi)=:
\sum_{I,J} f_{I,J}(z) \fra^J \xi^I
,$$
where $r_{I,J}$ are entire functions in $z$ and the summation is finite.
We define the map
$f\mapsto \pia( f)$ as above.


We define the {\it Fock--Berezin} space
$\SF_{p,q}(\cA)$
as the space of all  functions $f(z,\xi)$ satisfying
the condition: {\it for each $I$, $J$, the function $f_{I,J}(z)$
is in the Schwartz--Fock space $\cS\F_p$,} see Subsection
\ref{ss:schwartz-space}. We say that a sequence $f^{(k)}\in\SF_{p,q}(\cA)$ converges%
\footnote{This convergence corresponds to a topology of inductive limit.} to
$f$ if 

\sm

--- for all but a finite number of $J$ all $f^{(k)}_{I,J}$ are zero;

\sm

--- for each $I$, $J$, we have a convergence $f^{(k)}_{I,J}(z)\to f_{I,J}(z)$ in $\cS\F_p$.

\smallskip

{\sc Remark.}
 We   can assume that
all the $f_{I,J}$ are in the Hilbert--Fock space $\F_p$.
 But  Gauss--Berezin operators defined below 
can be unbounded in this space; therefore this point
 of view requires
descriptions of domains of
operators and examination of  products of operators.
 Our definition admits some variations
 (we chose an open dense subset in $\F_p$, and our  choice is volitional).
\hfill $\boxtimes$

\smallskip


{\bf\punct  Another form of the Gauss--Berezin integral.%
\label{ss:inner-products}}
 Consider
the integral
\begin{multline}
\int \exp\Bigl\{\frac 12
\begin{pmatrix}z&\xi\end{pmatrix}
\begin{pmatrix} A&B\\-B^t&C\end{pmatrix}
\begin{pmatrix}z^t\\ \xi^t\end{pmatrix} + z \alpha^t +\xi \beta^t \Bigr\}
\times   \\ \times
\exp\Bigl\{\frac 12
\begin{pmatrix}\ov z&\ov \xi\end{pmatrix}
\begin{pmatrix} K&L\\-L^t&M\end{pmatrix}
\begin{pmatrix}\ov z\\ \ov\xi\end{pmatrix} +
\ov z \kappa+\ov \xi \lambda^t
 \Bigr\}
\cdot
e^{-z\ov z^t -\xi\ov\xi^t}
\,dz\,\ov dz\,d\xi\,d\ov\xi
\label{eq:int-for-product}
,
\end{multline}
where two matrices $\begin{pmatrix} A&B\\-B^t&C\end{pmatrix}$, $\begin{pmatrix} K&L\\-L^t&M\end{pmatrix}$ 
have the same structure as in Subsect. \ref{ss:GB-integral}, row-vectors $\alpha$,
$\kappa$ are even, the vectors $\beta$, $\lambda$ are odd.

Since $\C^n\simeq\R^{2n}$,
this  integral is a special case of
the Gauss--Berezin integral.
We get
\begin{equation}
\mathrm{const}\cdot
\exp\left\{
\frac12
\begin{pmatrix} \alpha&\beta&\kappa&\lambda\end{pmatrix}
\begin{pmatrix}
-A&-B&1&0\\
B^t& -C&0&1\\
1&0& -K&-L\\
0&-1& L^t&-M
\end{pmatrix}^{-1}
\begin{pmatrix}
 \alpha^t\\ \beta^t\\ \kappa^t\\ \lambda^t
 \end{pmatrix}
\right\}
,
\label{eq:end-for-product}
\end {equation}
where the scalar factor is a hybrid of
the Pfaffian and  the Berezinian
mentioned above in Subsect. \ref{ss:evaluation-GB}.

\smallskip


{\bf\punct Integral operators.%
\label{ss:super-integral-operators}}
We write operators $\SF_{p,q}(\cA)\to\SF_{r,s}(\cA)$ as
\begin{equation}
A f(z,\xi)=
\int K(z,\xi; \ov u,\ov\eta)
\,f(u,\eta)\, e^{-z\ov z^t- \eta\ov\eta^t}
\,
du\,d\ov u\, d\ov\eta\,d\eta
\label{eq:integral-operator}
.
\end{equation}


{\bf\punct Linear and antilinear operators.%
\label{ss:antilinear-operators}}
We say that an operator $A:\SF_{p,q}(\cA)\to\SF_{r,s}(\cA)$
is {\it linear} if
$$
A(f_1+f_2)=Af_1+ Af_2,\qquad
A(\lambda f)=\lambda Af,\,\,\text{where $\lambda$ is a phantom constant}
,
$$
and {\it antilinear} if
$$
A(f_1+f_2)=Af_1+Af_2,\qquad
A(\lambda f)=\lambda^\sigma Af,\,\,\text{where $\lambda$ is a phantom constant}
,
$$
the automorphism $\lambda\mapsto\lambda^\sigma$.
Clearly, the operators
$$
Af(z,\xi)=\xi_j f(z,\xi),\qquad
B f(z,\xi)=\frac{\partial}{\partial\xi_j} f(z,\xi),\qquad
Cf(z,\xi)=\fra_j f(z,\xi)
$$
are antilinear.

\smallskip

An integral operator (\ref{eq:integral-operator}) is linear
if the kernel $K(z,\xi,\ov u, \ov\eta)$
is an even function in the total collection of all
Grassmann variables $\xi$, $\ov\eta$, $\fra$, i.e..
$$
K(z,\xi,\ov u,\eta)=K(z,-\xi,\ov u, -\ov\eta)^{\sigma}
.
$$
An operator (\ref{eq:integral-operator}) is antilinear if and only if the function $K$ is odd.

Below we meet only linear and antilinear operators.

\smallskip

We define also the (antilinear) operator $\mathsf S$ of $\sigma$-conjugation,
\begin{equation}
\mathsf S\,f(z,\xi) = f^\sigma(z,\xi).
\label{eq:sfS}
\end{equation}
Evidently,
$$
\mathsf S^2 f=f
.
$$

\smallskip


{\bf \punct Gauss--Berezin vectors in the narrow sense.%
\label{ss:vectors-narrow}}
A Gauss--Berezin vector (in the narrow sense)
is a vector of the form
\begin{equation}
\mathbf b \begin{bmatrix}
A&B\\ -B^t& C
\end{bmatrix}
=
\lambda
\exp\Bigl\{\frac 12
\begin{pmatrix} z&\xi\end{pmatrix}
\begin{pmatrix}
A&B\\ -B^t& C
\end{pmatrix}
\begin{pmatrix} z^t\\ \xi^t\end{pmatrix}
\Bigr\}
\label{eq:gauss-vector}
,
\end{equation}
where $A$, $B$, $C$ are  as above, see Subsection
\ref{ss:GB-integral}.

\begin{observation}
$\mathbf b[\cdot]\in \SF_{p,q}(\cA)$ if and only if
 $\|\pia(A)\|<1$.
\end{observation}


\smallskip


{\bf\punct Gauss--Berezin operators in the narrow sense.%
\label{ss:operator-narrow}}
A {\it  Gauss--Berezin integral operator in the narrow sense}
 is an integral operator
$$
\SF_{p,q}(\cA)\to\SF_{r,s}(\cA),
$$
whose kernel as a function in $z_1$, \dots, $z_r$, $\ov u_1$, \dots, $u_p$, $\xi_1$, \dots, $\xi_s$,
$\ov\eta_1$,
$\ov \eta_q$  is
 a Gauss--Berezin vector.
Precisely, a Gauss--Berezin operator
has the form
\begin{multline}
\B f(z,\xi)
=\\=
\lambda\cdot
\iint \exp\left\{
\frac 12
\begin{pmatrix}
z&\xi&\ov u&\ov \eta
\end{pmatrix}
\begin{pmatrix}
A_{11}&A_{12}&A_{13}&A_{14}\\
A_{21}&A_{22}&A_{23}&A_{24}\\
A_{31}&A_{32}&A_{33}&A_{34}\\
A_{41}&A_{42}&A_{43}&A_{44}
\end{pmatrix}
\begin{pmatrix}
z^t\\ \xi^t \\ \ov u^t \\ \ov\eta^t
\end{pmatrix}
\right\}\,f(u,\eta)
\times\\ \times
e^{-\eta\ov\eta^t-u\ov u^t}
\,du\,d\ov u\, d\eta\,
d\ov\eta
\label{eq:GBO-narrow}
,
\end{multline}
where $\lambda$ is an even phantom constant,
$A_{ij}$ is composed of even phantom constants if $(i+j)$ is even,
otherwise $A_{ij}$ is composed of odd phantom constants.
They also satisfy  the natural  symmetry conditions
for a matrix of a quadratic form in the variables
 $z$, $\xi$, $\ov u$, $\ov\eta$.

\smallskip

{\sc Remark.}
On the other hand, a Gauss--Berezin vector can be regarded
as a Gauss--Berezin operator
$\SF_{0,0}(\cA)\to\SF_{p,q}(\cA)$.
\hfill $\boxtimes$

\smallskip

For   Gauss--Berezin operators
$$
\B_1:\SF_{p,q}(\cA)\to \SF_{p',q'}(\cA),
\qquad \B_2:\SF_{p',q'}(\cA)\to \SF_{p'',q''}(\cA)
$$
evaluation of their product is reduced to
the Gauss--Berezin integral (\ref{eq:int-for-product}).
For operators  in  general position, we can apply
formula (\ref{eq:end-for-product}). Evidently, in this case the product is
a Gauss--Berezin operator again.
However, our final Theorem \ref{th:main-2}
avoids this calculation.

Also, considerations of  Section \ref{s:orthogonal}
suggest an extension of the definition
of Gauss--Berezin operators.

\smallskip


{\bf\punct General Gauss--Berezin operators.%
\label{ss:general-GB-operators}}
As above, we define  first order differential operators
$$
\frD[\xi_j]f:=
\Bigl(\xi_j+\frac \partial{\partial \xi_j}\Bigr)f
.
$$
If a function $f$ is independent of $\xi_j$, then
$$
\frD[\xi_j] f=\xi_j f, \qquad
\frD[\xi_j] \xi_j f= f
.
$$
Evidently,
$$
\frD[\xi_j]^2=1,\qquad
\frD[\xi_i]\,\frD[\xi_j]=-\frD[\xi_j]\,\frD[\xi_i],
\q, i\ne j
.
$$
The operators $\frD[\xi_j]$  are antilinear.

\smallskip

A {\it  Gauss--Berezin operator}
$\SF_{p,q}(\cA)\to\SF_{r,s}(\cA)$
is an operator of the form
\begin{equation}
\frC=\lambda
\frD[\xi_{i_1}]\dots \frD[\xi_{i_k}]
\,\B\,
\frD[\eta_{m_1}]\dots \frD[\eta_{m_l}]\cdot
\mathsf S^{k+l}
\label{eq:atlas}
,\end{equation}
where the operator $\mathsf S$
is given by (\ref{eq:sfS}) and

\smallskip

--- $\B$ is a Gauss --Berezin operator in the narrow sense;

\smallskip

---  $i_1<i_2<\dots<i_k$, $m_1<m_2<\dots<m_l$,
and $k$, $l\ge 0$;

\smallskip

---  $\lambda$ is an
even invertible phantom constant.

\smallskip

Note that a Gauss--Berezin operator is linear.

\smallskip



{\sc Remark.}  We define the set of Gauss--Berezin operators
as a  union of $2^{p+q}$ sets.
These sets are not disjoint.
Actually, we get a supermanifold consisting of two
connected components (according to the parity of $k+l$).
  Each set
(\ref{eq:atlas}) is open and dense in the corresponding component.
This will become obvious below.
\hfill $\boxtimes$


\smallskip

{\bf\punct Operators $\pia(\frB)$ and boundedness of Gauss--Berezin operators.%
\label{ss:remark}}
Let $K(z,\xi,\ov u,\ov\eta)$
be the kernel of a Gauss--Berezin operator.
Then the formula
$$
\frC f(z,\xi)=
\int \pia\Bigl(K(z,\xi; \ov u,\ov\eta)\Bigr)
\,f(u,\eta)\, e^{-z\ov z^t- \eta\ov\eta^t}
\,
du\,d\ov u\, d\ov\eta\,d\eta
$$
determines an integral operator
$$
\pia(\frC):\,\F_p\otimes\Lambda_q\to\F_r\otimes\Lambda_s
.
$$
Evidently, this operator is a tensor product of a Gaussian
operator
$$
\pia^+(\frC):\F_p\to\F_r
$$
 and a Berezin operator
$$
\pia^-(\frC):\Lambda_q\to\Lambda_s
.
$$
 For instance,
for an operator $\frB$ given by the standard formula (\ref{eq:GBO-narrow}),
we get  the Gaussian operator
$$
\pia^+(\frB)
f(z)=
\int_{\C^n}
\exp\Bigl\{ \frac 12
\begin{pmatrix}z&\ov u\end{pmatrix}
\biggl[\pia\begin{pmatrix}
A_{11}&A_{13}\\A_{31}&A_{33}
\end{pmatrix}\biggr]
\begin{pmatrix}z^t\\ \ov u^t\end{pmatrix}
\Bigr\}\,
f(z)\,e^{-z\ov z^t}\,dz\,d\ov z
$$
and the Berezin operator
$$
\pia^-(\frB)
g(\xi)=
\int
\exp\Bigl\{ \frac 12
\begin{pmatrix}\xi &\ov \eta\end{pmatrix}
\biggl[\pia\begin{pmatrix}
A_{22}&A_{24}\\A_{42}&A_{44}
\end{pmatrix}\biggr]
\begin{pmatrix}\xi^t\\ \ov \eta^t\end{pmatrix}
\Bigr\}\,
g(\xi)\,e^{-\xi\ov \xi^t}\,d\xi\,d\ov \xi
.$$

Next,
\begin{equation}
\pia^\pm(\frC_1\frC_2)=\pia^\pm(\frC_1)\,\pia^\pm(\frC_2)
\label{eq:otimes}
.
\end{equation}

\begin{observation}
 The Gauss--Berezin operator \eqref{eq:GBO-narrow} is bounded in the sense of Fock--Berezin spaces
 if and only if 
 the operator $\pi^+_\downarrow(\frB)$ is bounded {\rm(}i.e., satisfies conditions of Theorem {\rm\ref{th:olshanski})}. 
\end{observation}

{\sc Proof.} Only statement 'if' requires a proof.
A kernel of a Gauss--Berezin  operator has the form
$$
\sum_{I,J,K}\exp\Bigl\{\frac12 \begin{pmatrix}z&\ov u\end{pmatrix} 
\biggl[\pia
\begin{pmatrix}
A_{11}&A_{13}\\A_{31}&A_{33}
\end{pmatrix}\biggr]
\begin{pmatrix}z^t\\\ov u^t\end{pmatrix}
\Bigr\}
P_{I,J,K}(z,\ov u)\fra^I \xi^J \ov \eta^K,
$$
where $P_{I,J,K}$ is a polynomial in $z$, $\ov u$.
Denote the quadratic form in the exponential by $S(z)$ and decompose 
$P_{I,J,K}(z,\ov u)$ as a sum of monomials.
It is sufficient to show that the following operators  are bounded as operators $\cS \F_p\to\cS\F_r$:
$$
\int \prod z_j^{l_j} \exp\{S(z,\ov u)\} \prod \ov u_i^{k_i} f(u)\, e^{-z\ov z}\, d\lambda(u).
$$
Integrating by parts (see, e.g., \cite{Ner-book}, Subsect. V.3.5) we arrive at 
$$
\prod z_j^{l_j} \cdot \int \exp\{S(z,\ov u)\}\cdot \prod \Bigl( \frac\partial{\partial  u_i}\Bigr)^{k_i} f(u)\cdot e^{-z\ov z}\, d\lambda(u).
$$
It remains to notice that partial differentiations and multiplications by linear functions are bounded operators in the Fock-Schwartz spaces,
see Theorem \ref{th:gauss-schwartz}.
\hfill $\square$

\smallskip


{\bf\punct  Products of Gauss--Berezin operators.%
\label{ss:product-GB}}

\begin{theorem}
\label{th:main-0}
 For each  Gauss--Berezin operator
$$
\B_1:\SF_{p,q}(\cA)\to \SF_{p',q'}(\cA),
\qquad \B_2:\SF_{p',q'}(\cA)\to \SF_{p'',q''}(\cA)
,
$$
their product $\B_2 \B_1$
is
either a Gauss--Berezin operator
or
\begin{equation}
\pia(\B_2 \B_1  f)=0
\label{eq:product-caput}
\end{equation}
for all $f$.
\end{theorem}

Clearly, the condition \eqref{eq:product-caput} also is equivalent to $\pia^-(\B_2 \B_1  f)=0$.

\sm

For a proof, see
 Section \ref{s:correspondence}.
We also present an interpretation
of the product in terms of linear relations.

\smallskip

{\sc Remark.} In the case (\ref{eq:product-caput})
the kernel of the product has the form (\ref{eq:general})
but it is not a Gauss--Berezin operator in
our sense. Possibly, this requires
to change our   definitions.
\hfill$\boxtimes$

\smallskip


{\bf\punct General Gauss--Berezin vectors.}
A Gauss--Berezin vector is a vector
of the form
$$
\frD[\xi_{i_1}]\dots \frD[\xi_{i_k}]\,
\mathsf S^k
\b
,
$$
where $\b$ is a Gauss--Berezin vector in the narrow sense.

\smallskip

{\sc Remark.}
We write $\mathsf S^k$ as in 
(\ref{eq:atlas}). However, omitting this factor does not change the definition.
\hfill $\boxtimes$


\section{Supergroups $\OSp(2p|2q)$}

\label{s:group}

\COUNTERS

Here we define a super-analogue
of the groups $\OO(2n,\C)$ and $\Sp(2n,\C)$.
For a general exposition of supergroups
and super-Grassmannians, see books
\cite{Ber-super}, \cite{Man}, \cite{BL}, \cite{CF}.


\smallskip

{\bf\punct Modules $\cA^{p|q}$.%
\label{ss:superspaces}}
Let
$$
\cA^{p|q}:=\cA^p\oplus \cA^q
$$
be a direct sum of $(p+q)$ copies of $\cA$.
We regard elements  of
$\cA^{p|q}$
 as row-vectors
$$
(v_1,\dots,v_p; w_1,\dots,w_q)
.$$

We  define a structure of $\cA$-bimodule
on $\cA^{p|q}$.
The addition in $\cA^{p|q}$ is  natural.
The left multiplication by $\lambda\in\cA$
is also  natural
$$
\lambda\circ(v_1,\dots,v_p; w_1,\dots,w_q)
:=
(\lambda v_1,\dots,\lambda v_p;
\lambda w_1,\dots,\lambda w_q)
.
$$

The right multiplications by $\kappa\in\cA$ is
$$
(v_1,\dots,v_p; w_1,\dots,w_q)*\kappa
:=(v_1 \kappa,\dots,v_p\kappa;
 w_1\kappa^\sigma,\dots,w_q\kappa^\sigma)
,
$$
where $\sigma$ is the involution of $\cA$ defined above.

We define the {\it even part} of
$\cA^{p|q}$ as
$(\cA_{\even})^p\oplus(\cA_{\odd})^q$
and the {\it odd part} as
$(\cA_{\odd})^p\oplus(\cA_{\even})^q$.

\smallskip


{\bf \punct Matrices.%
\label{supermatrices}}
Denote by $\Mat(p|q;\cA)$ the space of
$(p+q)\times(p+q)$ matrices over $\cA$,
we represent such matrices in the block form
$$
Q=\begin{pmatrix} A&B\\ C&D\end{pmatrix}
.
$$
We say that a {\it matrix $Q$ is  even}  if
 all  matrix elements of $A$, $D$ are
even and all  matrix elements of $B$, $C$ are odd.
{\it A matrix is  odd} if elements of $A$, $D$ are odd
and elements of $B$, $C$ are even.

A matrix $Q$ acts on the space $\cA^{p|q}$
as
$$
v\to  v Q
.
$$
Such transformations   are compatible
with the  left $\cA$-module structure on $\cA^{p|q}$,
i.e.,
$$
(\lambda\circ v)\, Q= \lambda \circ (v Q)
\qquad \text{for any $\lambda\in\cA$,  $v\in\cA^{p|q}$}
.
$$
However, {\it even matrices also regard the right $\cA$-module
structure,}
$$
(v*\lambda)\, Q= (v Q)*\lambda
\qquad \text{for any $\lambda\in\cA$,  $v\in\cA^{p|q}$}.
$$
(we use the rule (\ref{eq:sigma!})).

\smallskip


{\bf\punct Super-transposition.%
\label{ss:super-transposition}}
The {\it supertranspose} of $Q$ is defined by
$$
Q^{st}
=\begin{pmatrix} A&B\\ C&D\end{pmatrix}^{st}
:=\begin{cases}
\begin{pmatrix}
A^t& C^t\\-B^t& D^t
\end{pmatrix} \quad \text{if $Q$ is even}
,
\\
\begin{pmatrix}
A^t& -C^t\\B^t& D^t
\end{pmatrix} \quad \text{if $Q$ is odd}
,
\end{cases}
$$
and
$$
(Q_1+Q_2)^{st}
:=
Q_1^{st}+Q_2^{st}.
$$

The following identity holds
\begin{equation}
(Q R)^{st}=\begin{cases}
R^{st}Q^{st}\quad & \text{if $Q$ or $R$ are even}
,
\\
-R^{st}Q^{st}\quad & \text{if both $R$ and $Q$ are odd}
.
\end{cases}
\label{eq:product-transposition}
\end{equation}

Below we  use only the first row.

\smallskip


{\bf\punct The supergroups $\GL(p|q;\cA)$.
\label{ss:supergroup}}
The group $\GL(p|q;\cA)$ is the group of
{\it even} invertible
matrices in $\Mat(p|q;\cA)$. The following lemma
is  trivial.

\begin{lemma}
\label{l:super-invertible}
 An even matrix
$Q\in\Mat(p|q;\cA)$ is invertible

\smallskip

{\rm a)} if and only if the matrices $A$, $D$
are invertible;

\smallskip

{\rm b)}  if and only if
the matrices $\pia(A)$, $\pia(D)$ are invertible.
\end{lemma}

Here $\pia(A)$ denotes the matrix composed of elements
$\pia(a_{kl})$.

Also, the map $Q\mapsto \pia(Q)$ is a well-defined epimorphism
$$
\pia: \GL(p|q;\cA)\to\GL(p,\C)\times\GL(q,\C)
$$
(because $\pia(B)=0$, $\pia(D)=0$).

\smallskip


{\bf\punct The supergroup $\OSp(2p|2q;\cA)$.%
\label{ss:OSp}} We define the standard
{\it orthosymplectic form} on $\cA^{2p|2q}$
by
$$
\fros(u,v):=u J v^{st}
,
$$
where $J$ is a block $(p+p+q+q)\times (p+p+q+q)$
matrix
\begin{equation}
J:=
\frac12
\begin{pmatrix}0&1_p&0&0\\
              -1_p&0&0&0\\
               0&0&0&1_q\\
               0&0&1_q&0
\end{pmatrix}
.
\label{eq:J}
\end{equation}
The group $\OSp(2p|2q;\cA)$ is
the subgroup in $\GL(2p|2q;\cA)$ consisting of matrices
$g$
satisfying
$$
\fros(u,v)=\fros(ug,vg)
.
$$
Equivalently,
$$
g J g^{st}=J
$$
(for this conclusion,
we use (\ref{eq:product-transposition}); since $g\in \GL(p|q;\cA)$
is even,
a sign does not appear).

We also write elements of $\OSp(2p|2q;\cA)$ as block $(2p+2q)$-matrices
$g=\begin{pmatrix}A&B\\C&D\end{pmatrix}$. For such matrix, we have
\begin{align*}
\pia(A)\begin{pmatrix}0&1_p\\-1_p&0\end{pmatrix}\pia(A)^t&=
\begin{pmatrix}0&1_p\\-1_p&0\end{pmatrix},\qquad\text{i.e., $\pia(A)\in \Sp(2n,\C)$};
\\
\pia(D) \begin{pmatrix}0&1_q\\1_q&0\end{pmatrix}\pia(D)^t&=\begin{pmatrix}0&1_q\\1_q&0\end{pmatrix},
\qquad\text{i.e., $\pia(D)\in \OO(2n,\C)$}
.
\end{align*}


{\bf\punct The super-Olshanski
 semigroup $\Gamma\OSp(2p|2q;\cA)$.%
\label{ss:olshanski-super}}
We define the  semigroup 
$\Gamma\OSp(2p|2q;\cA)$ as a subsemigroup
in $\OSp(2p|2q;\cA)$ consisting of matrices
$g=\begin{pmatrix} A&B\\C&D\end{pmatrix}$
such that
$\pia(A)$ is contained in the Olshanski
semigroup $\Gamma\Sp(2p,\R)$, see
Subsection \ref{ss:olshanski}.



\section{Super-Grassmannians}

\label{s:grassmannian}

\COUNTERS

This section is a preparation
to the definition of super-linear relations.


\smallskip

{\bf\punct Super-Grassmannians.%
\label{ss:super-grassmannian}}
 Let $u_1$,\dots, $u_r$
be even vectors and $v_1$,\dots, $v_s$
be odd vectors in $\cA^{p|q}$.
We suppose that

\smallskip

--- $\pia(u_j)\in(\C^p\oplus 0)$ are linearly independent,

\smallskip

--- $\pia(v_k)\in(0\oplus \C^q)$ are linearly independent.

\smallskip


{\it A supersubspace of superdimension $r|s$} is
a left $\cA$-module generated by such vectors.
Subspaces also are right $\cA$-submodules.

\smallskip

We define the {\it super-Grassmannian} $\Gr^{r|s}_{p|q}(\cA)$
as the space of all  supersubspaces in
$\cA^{p|q}$ of superdimension
$r|s$.

\smallskip

By the definition, the map $\pia$ projects
$\Gr^{r|s}_{p|q}(\cA)$ to the product
$\Gr^r_p\times \Gr^s_q$ of
the usual complex Grassmannians.
We denote by $\pia^\pm$ the natural projections
$$
\pia^+:\Gr^{r|s}_{p|q}(\cA) \to
\Gr_p^r,\qquad \pia^-:\Gr^{r|s}_{p|q}(\cA)\to\Gr^s_q
.
$$



{\bf\punct Intersections of subspaces.%
\label{ss:intersections}} Let us examine superdimensions of intersections of subspaces.

\begin{lemma}
\label{l:intersection}
Let $L$ be a subspace of superdimension $r|s$ 
in $\cA^{p|q}$, $M$ be a subspace of superdimension $\rho|\sigma$. Let
the following transversality conditions hold
\begin{equation}
\pia^+(L) + \pia^+(M)=\C^p,
\q
\pia^-(L) + \pia^-(M)=\C^q
\label{eq:transversality}
.
\end{equation}
Then
$L\cap M$ is a subspace and its superdimension
is
$(r+\rho-p)|(s+\sigma-q)$.
\end{lemma}

{\sc Remark.} If the transversality
conditions are not satisfied,
then incidentally  $L\cap M$
is not a subspace.
For instance, consider $\cA^{1|1}$
with a basis $e_1$, $e_2$ and subspaces
$$
L:=\cA(e_1+\fra_1 e_2),\qquad M:=\cA\cdot e_1
.
$$
Then $L\cap M=\cA \fra_1 e_1$ is not a subspace.
 \hfill $\boxtimes$

\smallskip

{\sc Proof.} Denote by $I\subset\cA$ the ideal
spanned by all $\fra_j$, i.e., $\cA/I=\C$.

It is easy to see that
\begin{equation}
L+M=\cA^{p|q}
\label{eq:L+M}
.
\end{equation}
Indeed, denote by $e_1$, \dots, $e_p$, $e_{p+1}$, \dots, $e_{p+q}$
the standard basis in $\C^p\oplus \C^q$. Then for each $k$ the submodule
 $L+M$ contains a vector
of the form
$$
E_k=e_k +\sum_{J\ne \emptyset} \sum_m x_{k,m,J} \fra^J e_m,
$$
where $x_{k,m,J}\in \C$. Actually, only a finite number of nonzero constants $\fra_j$
are contained in this expression. Without loss of generality, we can assume 
that this set is $\fra_1$, \dots, $\fra_N$. Then $L+M$ contains all vectors
$\fra_1\dots\fra_N E_k=\fra_1\dots\fra_N e_k$ and therefore
$L+M\supset  \fra_1\dots\fra_N\cdot \C^p\oplus\C^q$. Next, for each $l$
$$
\fra_1\dots \fra_{l-1}\fra_{l+1}\dots   \fra_N\, E_k
- \fra_1\dots \fra_{l-1}\fra_{l+1} \dots  \fra_N \,e_k
\,\,\in \,\, \fra_1\dots\fra_N \cdot\C^p\oplus\C^q.
$$
Therefore, $\fra_1\dots \fra_{l-1}\fra_{l+1} \dots   \fra_N \cdot\C^p\oplus\C^q\subset L+M$. Repeating this process, we get $\C^p\oplus \C^q\subset L+M$ and this implies \eqref{eq:L+M}.

Let $v\in\pia^+(L)\cap  \pia^+(M)$.
Choose
$x\in L$, $y\in M$ such that
$\pia(x)=v$, $\pia(y)=v$.
Then $x-y\in I\cdot \cA^{p|q}$.
However,
$I\cdot L+I\cdot M=I\cdot \cA^{p|q}$,
therefore we can represent
$$
x-y=a-b,\qquad\text{where $a\in I\cdot L$, $b\in I\cdot M$}
.
$$
Then
$$
(x-a)\in L, \,(y-b)\in M,\,\,
\pia(x-a)=v=\pia(y-b)
.
$$
Thus, for any vector $v\in\pia^+(L)\cap  \pia^+(M)$,
there is a vector $v^*\in L\cap M$ such that
$\pia(v^*)=v$.
The same is valid for vectors
$w\in\pia^-(L)\cap  \pia^-(M)$.

\smallskip

Therefore, $L\cap M$ contains a supersubspace of desired superdimension
generated by vectors $v^*$, $w^*$.
It remains to show that there are no extra vectors in the intersection.

\smallskip

Now, let us vary a phantom algebra $\cA$.
If $\cA$ is an algebra in a finite number of Grassmann
constants $\fra_1$,\dots, $\fra_n$, then
this
completes a proof, since (\ref{eq:L+M})
gives the same superdimension of the intersection
 over $\C$.

\smallskip

Otherwise, we choose a basis in $L$ and a basis in $M$.
Expressions for basis vectors contain only a finite
number of Grassmann constants $\fra_1$, \dots, $\fra_k$.
After this, we apply the same reasoning to
algebras $\cA[l]$ generated by Grassmann constants
$\fra_1$, \dots, $\fra_l$ for all $l\ge k$ and
observe that $L\cap M$ does not contain extra vectors.
\hfill $\square$

%
%
%
%
%

\smallskip


{\bf\punct Atlas on the super-Grassmannian.%
\label{ss:atlas-super}}
Define an atlas on the super-Grassmannian
$\Gr^{r|s}_{p|q}(\cA)$ as usual. Namely,
consider the following complementary subspaces
$$
V_+:=(\cA^r\oplus 0)\oplus (\cA^s\oplus 0)
\q\q
V_-:=(0\oplus\cA^{p-r})\oplus (0\oplus\cA^{q-s})
$$
in $\cA^{p|q}$. Let $S:V_+\to V_-$ be
an even operator.
 Then its graph
 is an element
of the super-Grassmannian.

Permuting coordinates in $\cA^p$ and $\cA^q$,
we get an atlas that covers the whole  super-Grassmannian $\Gr^{r|s}_{p|q}(\cA)$.

\smallskip


{\bf\punct Lagrangian super-Grassmannians.%
\label{ss:super-lagrangian}}
Now, equip the space $\cA^{2p|2q}$ with the
orthosymplectic form $\fros$ as above.
We say that a subspace $L$ is {\it isotropic}
if the form $\fros$ is zero on $L$.
A {\it Lagrangian subspace} $L$ is an isotropic
subspace of the maximal possible superdimension, i.e., $\dim L=p|q$.

\begin{observation}
Let $L$ be a super-Lagrangian subspace.
Then

\smallskip

 --- $\pia^+(L)$ is
Lagrangian
subspace in $\C^{2p}$ with respect to the
skew-symmetric bilinear form
$\begin{pmatrix} 0&1\\-1&0 \end{pmatrix}$;

\smallskip

--- $\pia^-(L)$ is
a Lagrangian subspace in $\C^{2q}$ with respect
to the symmetric bilinear form
$\begin{pmatrix} 0&1\\1&0 \end{pmatrix}$.
\end{observation}

\smallskip


{\bf\punct Coordinates on Lagrangian super-Grassmannian.}
Consider the following complementary Lagrangian subspaces
\begin{equation}
V_+:=(\cA^p\oplus 0)\oplus (\cA^q\oplus 0)
,
\q\q
V_-:=(0\oplus\cA^{p})\oplus (0\oplus\cA^{q})
\label{eq:v+v-}
.
\end{equation}

\begin{proposition}
\label{pr:karta-for-super}
Consider an even operator
$S:V_+\to V_-$,
$$
S=\begin{pmatrix} A&B\\C&D\end{pmatrix}
.
$$
The graph of $S$ is a Lagrangian subspace
if and only if
\begin{equation}
A=A^t,\q D=-D^t, \q C+ B^t=0
\label{eq:conditions}
.
\end{equation}
\end{proposition}

{\sc  Remark.} This statement is a super-imitation
of Lemma \ref{l:potapov}.
\hfill $\boxtimes$

\smallskip

{\sc Proof.} We write out
a vector $h\in\cA^{2p|2q}$ as
$$h=(u_+,u_-;v_+,v_-)\in
\cA^p\oplus\cA^p\oplus\cA^q\oplus\cA^q
.$$
Then
$$
\fros(h',h)=
u_+'\, (u_-)^{st}- u_-'\, (u_+)^{st}
+
v_+'\, (v_-)^{st}+ v_-'\, (v_+)^{st}
.$$
Let $h$ be in the graph of $S$.
Then
$$
\begin{pmatrix}
u_-& v_-
\end{pmatrix}=
\begin{pmatrix}
u_+ & v_+
\end{pmatrix}
\begin{pmatrix} A&B\\C&D\end{pmatrix}
=
\begin{pmatrix}
u_+ A +  v_+ C
&
u_+ B +  v_+ D
\end{pmatrix}
$$
and
\begin{multline*}
\fros(h',h)=
u_+' (u_+ A +  v_+ C)^{st}
-
(u_+' A +  v_+' C) u_+^{st}
+\\+
v_+'(u_+ B +  v_+ D)^{st}
+
(u'_+ B +  v'_+ D) v_+^{st}
.
\end{multline*}

Observe
 that the matrices $A$, $B$, $C$, $D$ are even%
\footnote{Recall that this means that $A$, $D$
are  composed of even phantom constants
and $C$, $B$ of odd phantom constants.}%
;
for this reason, we write $(u_+ A)^{st}=A^{st}(u_+)^{st}$
etc.,
see (\ref{eq:product-transposition}).
We arrive at
\begin{multline*}
u_+' \bigl[A^{st} (u_+)^{st}  +  C^{st} (v_+)^{st} \bigr]
-
\big[(u_+' A +  v_+' C\bigr] u_+^{st}
+\\+
v_+'\bigl[B^{st} (u_+)^{st}  +  D^{st} (v_+)^{st} \bigr]
+
\bigl[u'_+ B +  v'_+ D\bigr] v_+^{st}
.
\end{multline*}
Next,
$$
A^{st}=A^t,\q B^{st}=-B^t,\q C^{st}= C^t,\q D^{st}=D^t
.
$$
Therefore, we convert our expression to
the form
$$
u_+'(A-A^t) (u_+)^{st}
+
v_+'(D+D^t) (v_+)^{st}
+
(u_+')(B^t+C)v_+^{st}
+
(v^+) (B+C^t) u_+^{st}
.$$
This expression is zero if and only if the conditions
 (\ref{eq:conditions}) are satisfied.
\hfill $\square$

\smallskip


{\bf\punct Atlas on the Lagrangian super-Grassmannian%
\label{ss:latlas}.}
Now we imitate the construction of
Subsection
\ref{ss:atlas-lagrange0}.

Consider the standard basis in $\cA^{2p|2q}$ consisting 
of vectors, whose coordinates are 0 except one unit.
Denote elements of this basis by
$$
e_1,\dots,e_p;\, e'_1,\dots,e'_p;\, f_1,\dots,f_q;\, f'_1,\dots,f'_q.
$$
In this basis, the matrix of the orthosymplectic form is \eqref{eq:J}.
Consider   subsets
$I\subset\{1,2,\dots,p\}$,
$J\subset \{1,2,\dots,q\}$.

 We define
\begin{align}
V_+[I,J]=\bigl(\oplus_{i\in I} \cA e_i\bigr)
\oplus
\bigl(\oplus_{k\notin I} \cA e_k'\bigr)
\oplus
\bigl(\oplus_{j\in J} \cA f_j\bigr)
\oplus
\bigl(\oplus_{l\notin J} \cA f_l'\bigr)
\label{eq:la-}
,
\\
V_-[I,J]=\bigl(\oplus_{i\notin I} \cA e_i\bigr)
\oplus
\bigl(\oplus_{k\in I} \cA e_k'\bigr)
\oplus
\bigl(\oplus_{j\notin J} \cA f_j\bigr)
\oplus
\bigl(\oplus_{l\in J} \cA f_l'\bigr)
\label{eq:la+}
.
\end{align}
We denote by $\cO[I,J]$ the set of
all  the Lagrangian subspaces
that are graphs of even operators
$$
S:V_+[I,J]\to V_-[I,J]
.
$$
In fact, these operators satisfy the same conditions
as in Proposition  \ref{pr:karta-for-super}
(our initial chart is $\cO[\emptyset,\emptyset]$).
Thus, we get an atlas on the
Lagrangian super-Grassmannian.

\smallskip


{\bf\punct Elementary reflections.%
\label{ss:elementary-2}}
Now we repeat considerations of Subsection \ref{ss:elementary}.
We define elementary reflections
$\sigma[e_i]$, $\sigma[f_j]$ in $\cA^{2p|2q}$
by
\begin{align*}
\sigma[e_i]\, e_i^+= - e_i^-,\quad \sigma[e_i]\, e_i^-=e_i^+,\quad
\\
\sigma[e_i]\, e_k^\pm=e_k^\pm\,\, \,\,\text{for $k\ne i$},
\quad \sigma[e_i]\, f_j^\pm=f_j^\pm
,\end{align*}
and
\begin{align}
\sigma[f_j]\, f_j^+=  f_j^-,\quad \sigma[f_j]\, f_j^-=f_j^+,\quad
\label{eq:sigma-1}
\\
\sigma[f_i]\, f_k^\pm=e_k^\pm\,\,\,\,
 \text{for $k\ne j$},\quad \sigma[f_j]\, e_i^\pm=e_i^\pm
\label{eq:sigma-2}
,\end{align}
in the first row, we have an extra change of a sign because we want to
preserve the symplectic form.
Then
$$
\cO[I,J]=
\prod_{i\in I}\sigma[e_i]\cdot
\prod_{j\in J}\sigma[f_j]\cdot
\cO[\emptyset,\emptyset ]
.$$


\section{Superlinear relations}

\label{s:relations}

\COUNTERS

Gauss--Berezin integral operators are
enumerated by contractive Lagrangian super-linear relations.
These objects are defined in this section.

\smallskip


{\bf \punct  Superlinear relations.%
\label{ss:super-relations}}
We define  super-linear relations
$P:\cA^{p|q}\arr\cA^{r|s}$
as subspaces in $\cA^{p|q}\oplus\cA^{r|s}$.
Products are defined  as above, see
Subsection \ref{ss:prl}.

Next, for a superlinear relation we define
complex linear relations
$$
\pia^+(P):\C^p\arr\C^r,
\q
\pia^-(P):\C^q\arr\C^s
$$
in the natural way, we simply project
the super-Grassmannian in $\cA^{p|q}\oplus \cA^{r|s}$
onto the product of the complex Grassmannians.

\smallskip


{\bf\punct Transversality conditions.
\label{ss:transversality}}
Let $V$, $W$, $Y$ be {\it complex} linear spaces. We say that
linear relations
$$
P:V\arr W,\,\,Q:W\arr Y
$$
are {\it transversal} if
\begin{align}
\im P+\dom Q=W
\label{eq:transversality-1}
,
\\
 \indef P\cap \ker Q=0
 \label{eq:transversality-2}
.
\end{align}

We met these conditions  in Section \ref{s:orthogonal}, in what follows they are even more
important.

\begin{theorem}
\label{th:dim-product-1}
If $P:V\arr W$, $Q:W\arr Y$ are transversal,
then
$$
\dim QP = \dim Q+\dim P-\dim W
.
$$
\end{theorem}

{\sc Proof.} We rephrase the definition
of the product $QP$ as follows (see \cite{Ner-book}, Prop. II.7.1).
Consider the space $V\oplus W\oplus W\oplus Y$
and the following
 subspaces

\smallskip

--- $P\oplus Q$,

\smallskip

--- the subspace $H$ consisting of vectors
    $v\oplus w\oplus w\oplus y$,

\smallskip

--- the subspace $T\subset H$ consisting of
vectors $0\oplus w\oplus w\oplus 0$.

\smallskip

Let us project $(P\oplus Q)\cap H$
 on $V\oplus W$ along $T$.
The result is $QP\subset V\oplus W$.

By the first transversality condition (\ref{eq:transversality-1}),
$$
(P\oplus Q)+ H=V\oplus W\oplus W\oplus Y
,$$
therefore we know the superdimension
of the intersection $S:=(P\oplus Q)\cap H$.

\smallskip

By the second condition (\ref{eq:transversality-2})
the projection
$H\to V\oplus W$ is injective on $S$.
\hfill
$\square$

\smallskip


{\bf\punct Transversality for super-linear relations.%
\label{ss:super-transversality}}
We say that super-linear relations $P:V\arr W$
and $Q:W\arr Y$ are {\it transversal} if
$\pia^{+}(P)$ is transversal to $\pia^{+}(Q)$
and
$\pia^{-}(P)$ is transversal to $\pia^{-}(Q)$.

\begin{theorem}
\label{th:dim-product-2}
If
 $P:V\arr W$, $Q:W\arr Y$ are transversal
 super-linear relations,
then their product is a super-linear relation and
$$
\dim QP
 = \dim Q+\dim P-\dim W
.
$$
\end{theorem}

{\sc Proof.} We  follow
the proof of the previous theorem.
\hfill $\square$

\smallskip


{\bf\punct Lagrangian super-linear relations.%
\label{ss:lagrangian-super-relations}}
Consider the spaces $V=\cA^{2p|2q}$, $W=\cA^{2r|2s}$
endowed with the orthosymplectic forms
$\fros_V$, $\fros_W$, respectively.
Define the form
$\fros^\ominus$ on $V\oplus W$ as
$$
\fros^\ominus(v\oplus w,v'\oplus w'):=
\fros_V(v,v')-\fros_W(w,w')
.$$
A Lagrangian super-linear relation $P:V\arr W$
is a Lagrangian
supersubspace in $V\oplus W$.

\begin{observation}
\label{obs:graph-osp}
Let $g\in \OSp(2p|2q;\cA)$. Then the graph
of $g$ is a Lagrangian super-linear relation
$\cA^{2p|2q}\arr\cA^{2p|2q}$.
\end{observation}

\begin{theorem}
\label{ss:osp-category}
Let $P:V\arr W$, $Q:W\arr Y$ be
transversal Lagrangian super-linear
relations. Then
$QP:V\arr Y$ is  a Lagrangian
super-linear relation.
\end{theorem}

{\sc Proof.} Let $v\oplus w$, $v'\oplus w'\in P$
and $w\oplus y$, $w'\oplus y'\in Q$.
By definition,
$$
\fros_V(v,v')=\fros_W(w,w')=\fros_Y(y,y')
,
$$
therefore $QP$ is isotropic.
By the virtue of Theorem
\ref{th:dim-product-2}, we know $\dim QP$.
\hfill $\square$

\smallskip


{\bf\punct Components of Lagrangian super-Grassmannian.%
\label{ss:super-components}}
As we observed in Subsection \ref{ss:components-orthogonal},
the orthogonal Lagrangian Grassmannian in the space $\C^{2n}$
consists of two components. The usual
symplectic Lagrangian Grassmannian
is connected.
Therefore, the Lagrangian super-Grassmannian
consists of two components.

Below we must distinguish them.

Decompose $V=V_+\oplus V_-$, $W=W_+\oplus W_-$ as above
(\ref{eq:v+v-}). We say that the component
containing the linear relation
$$
(V_+\oplus W_-): \,V\arr W
$$
is even; the other component is odd.


\smallskip

{\bf\punct Contractive Lagrangian linear relations.%
\label{ss:super-contractive}}
Now, we again (see Section \ref{s:symplectic}) consider
the Hermitian form $M$ on $\C^{2p}$, it is defined
by a matrix
$\begin{pmatrix}1_p&0\\0&-1_p \end{pmatrix}$.
Then $\C^{2p}$ becomes an object of the category
$\bSp$.

We say that a Lagrangian super-linear relation
$P:V\arr W$
is  contractive if  $\pia^+(P)$ is a morphism
of the category $\bSp$.


\smallskip

{\bf\punct Positive domain in the Lagrangian super-Grassmannian.%
\label{ss:positive-domain}}
We say that a Lagrangian subspace $P$ in $\C^{2p|2q}$
is positive if the form $M$ defined in the prevous subsection is positive
on $\pia^+(P)$.


\section{Correspondence between\\ Lagrangian superlinear
relations\\ and Gauss--Berezin operators}

\label{s:correspondence}

\COUNTERS

Here we prove our main results, namely Theorems \ref{th:main-1},
\ref{th:main-2}.


\smallskip

{\bf\punct Creation--annihilation operators.%
\label{ss:super-creation}}
Let $V:=\cA^{2p|2q}$ be a superlinear space
endowed with the  orthosymplectic
bilinear form $\frs$ defined by the matrix \eqref{eq:J}. For a vector
$$
v\oplus w:=v_+\oplus v_-\oplus w_+\oplus w_-\in \cA^{2p|2q}
,
$$
we define the creation-annihilation
operator in the Fock--Berezin space
$\SF_{p,q}(\cA)$ by
$$
\wh a(v\oplus w)f(z,\xi)=\Bigl(\sum_iv_+^{(i)}
\frac {\partial}{\partial z_i}+
\sum_i v_-^{(i)} z_i
+
\sum_j w_+^{(j)}
\frac {\partial}{\partial \xi_i}
+
\sum_j  w_-^{(j)} \xi_j
\Bigr)\, f(z,\xi)
.
$$


{\bf\punct Supercommutator.%
\label{ss:super-commutator}}
We say that a vector $v\oplus w$ is even if $v$ is even and $w$ is odd.
It is odd if $v$ is odd and $w$ is even. This corresponds
to the definition of even/odd  for
$(1|0)\times(2p|2q)$ matrices.
Let $h=v\oplus w$, $h'=v'\oplus w'$.
We define the supercommutator
$[\wh a (h),\wh a(h')]_{\super}$
as
$$
[\wh a (h),\wh a(h')]_{\super}
=
\begin{cases}
[\wh a (h),\wh a(h')]=\wh a (h)\wh a(h')-\wh a (h')\wh a(h)
\quad\text{if  $h$ or $h'$ is even,}
\\
\{\wh a (h),\wh a(h')\}=\wh a (h)\wh a(h')+\wh a (h')\wh a(h)
\quad\text{if  $h$, $h'$ are odd}
.
\end{cases}
$$

Then
$$
[\wh a (h),\wh a(h')]_{\super}
=
\frs(h,h')\cdot 1
,$$
where $1$ denotes the unit operator.

\smallskip

Also, note that an operator $\wh a(h)$ is linear
(see Subsection \ref{ss:antilinear-operators})
if $h$ is even and antilinear if $h$ is odd.

\smallskip


{\bf\punct Annihilators of Gaussian vectors.%
\label{ss:anigilyatory}}

\begin{theorem}
\label{th:vector-annihilator}
{\rm a)} For a  Gauss--Berezin vector
$\mathbf b\in \SF_{p,q}(\cA)$ consider the set $L$ of all
 vectors $h\in\cA^{2p|2q}$ such that
$$
\wh a(h)\, \mathbf b=0
.
$$
Then $L$ is a positive Lagrangian subspace
in $\cA^{2p|2q}$.

\smallskip

{\rm b)} Moreover, the map $\mathbf b\mapsto L$
is a bijection
$$
\left\{
\begin {matrix}
\text{The set of
all  Gauss--Berezin vectors}\\
\text{defined up to an invertible scalar}
\end{matrix}
\right\}
\leftrightarrow
\left\{\begin{matrix}
 \text{The positive}
\\
\text{ Lagrangian Grassmannian}
\end{matrix}
\right\}
.
$$
\end{theorem}


\smallskip

Before we begin a formal proof we propose the following
(insufficient, but clarifying) argument.
Let $h$, $h'\in L$. If one of them is even, then
we write
$$
\Bigl(\wh a (h)\wh a(h')-\wh a (h')\wh a(h)\Bigr)\, \mathbf b
.
$$
By the definition of $L$, this is 0. On the other hand,
this is $\frs(h,h')\mathbf b$. Therefore, $\frs(h,h')=0$.

If both $h$, $h'$ are odd, then we write
$$
0=\Bigl(\wh a (h)\wh a(h')+\wh a (h')\wh a(h)\Bigr)\,\mathbf b
=\frs(h,h')\b
$$
and arrive at the same result.

\smallskip

{\sc Proof.} First, let $\b(z,\xi)$
have the standard form (\ref{eq:gauss-vector}).
We write out
\begin{multline*}
\wh a(v\oplus w) \b(z,\xi)=
\Bigl(\sum_iv_+^{(i)}
\frac {\partial}{\partial z_i}+
\sum_i v_-^{(i)} z_i
+
\sum_j w_+^{(j)}
\frac {\partial}{\partial \xi_i}
+
\sum_j  w_-^{(j)} \xi_j
\Bigr)
\times\\ \times
\exp\Bigl\{\frac 12
\begin{pmatrix} z&\xi\end{pmatrix}
\begin{pmatrix}
A&B\\ -B^t& C
\end{pmatrix}
\begin{pmatrix} z^t\\ \xi^t\end{pmatrix}
\Bigr\}
=\\=
\Bigl( v_+ (Ax^t+B\xi^t)+ v_- z^t+
w_+ (-B^t z^t+C\xi^t)+ w_- \xi^t
\Bigr)\cdot \b(z,\xi)
=\\=
\Bigl( (v_+ A- w_+ B^t+ v_-)z^t+
(v_+ B+w_+ D + w_-) \xi^t\Bigr) \cdot \b(z,\xi)
.
\end{multline*}

This is zero if and only if
$$
\begin{cases}
v_-=-(v_+ A- w_+ B^t)\\
 w_-=-( v_+ B+w_+ D)
\end{cases}
.
$$
However, this system of equations determines
a Lagrangian subspace.
The positivity of a Lagrangian subspace is equivalent
to $\|\pia(A)\|<1$ (see, for instance \cite{Ner-lectures}).

\smallskip

Next, consider an arbitrary Gauss--Berezin vector
\begin{equation}
\mathbf b(z,\xi)
=
\frD[\xi_{i_1}]\dots \frD[\xi_{i_k}]
\,\mathsf S^k\,\,
\b[T]
\label{eq:42}
,
\end{equation}
where $\b[T]$
is a standard Gauss--Berezin vector.
We have
\begin{equation}
\wh a\bigl(h)\,
\frD[\xi_1]\, \mathsf S=
\frD[\xi_1]\, \mathsf S\,
\wh a\bigl(\sigma[f_1]h\bigr)
\label{eq:43}
,
\end{equation}
where
$\sigma[f_1]$
is an elementary reflection
given by (\ref{eq:sigma-1})--(\ref{eq:sigma-2}).

If
$h$
ranges in a Lagrangian subspace, then
$\sigma[f_1]h$
also ranges in (another) Lagrangian subspace.
Also, a map $\sigma[f_1]$ takes positive subspaces
to positive subspaces.
Therefore, the statements a) for vectors
$$
\frD[\xi_{i_1}]\, \frD[\xi_{i_2}] \dots \frD[\xi_{i_k}]
\,\mathsf S^k\,\,
\b[T]\qquad\text{and}\qquad
\frD[\xi_{i_2}]\dots \frD[\xi_{i_k}]
\,\mathsf S^{k-1}\,\,
\b[T]
.$$
are equivalent.

In fact, for fixed $i_1$, \dots,  $i_k$,
all  vectors of the form (\ref{eq:42})
correspond to a fixed chart in the Lagrangian
super-Grassmannian,
namely to
$$
 \sigma[f_{i_1}]\cdots \sigma[f_{i_k}]
\cdot
\cO[\emptyset,\emptyset]
$$
in notation of Subsection \ref{ss:elementary-2}.

But these charts cover the set of all positive Lagrangian
subspaces.
\hfill $\square$

\smallskip


\begin{theorem}
\label{th:uniqueness}
For a  positive Lagrangian subspace
$L\subset \cA^{2p|2q}$,
consider the system
of equations
\begin{equation}
\wh a(v\oplus w) f(z,\xi)=0
\q\text{for all $v\oplus w\in L$}
\label{eq:41}
,
\end{equation}
for a function $f(z,\xi)$.
All its solutions are of the form
$\lambda\b(z,\xi)$, where $\b(z,\xi)$
is a Gauss--Berezin vector and $\lambda$ is a phantom constant.
\end{theorem}

{\sc Proof.} It suffices to  prove the statement
for $L$  in the principal chart.
Put
$$
\phi(z,\xi):= f(z,\xi)/\b(z,\xi)
,
$$
i.e.,
$$
f(z,\xi)=\b(z,\xi)\cdot \phi(z,\xi)
.$$
 By the Leibniz rule,
\begin{multline*}
0=
\wh a(v\oplus w)
\bigl(
\b(z,\xi) \phi(z,\xi)
\bigr)
=\\=
\Bigl(
\wh a(v\oplus w)\b(z,\xi)\Bigr)
\cdot \phi(x,\xi)
+
\b(z,\xi)\cdot
\Bigl(
\sum_j v_+^{(j)}
\frac {\partial}{\partial z_i}
+
\sum_j w_+^{(j)}
\frac {\partial}{\partial \xi_i}
\Bigr)
.\phi(z,\xi)
\end{multline*}

The first summand is zero by the definition
of $\b(z,\xi)$. Since $v_+$, $w_+$ are arbitrary,
we get
$$
\frac {\partial}{\partial z_i} \phi(z,\xi)=0,\q
\frac {\partial}{\partial \xi_i} \phi(z,\xi)=0
.
$$
Therefore, $\phi(z,\xi)$ is a phantom constant.
\hfill$\square$

\smallskip


{\bf\punct Gauss--Berezin operators and
superlinear relations.%
\label{ss:GB-relations}}
Let $V=\cA^{2p|2q}$, $\wt V=\cA^{2r|2s}$
be (super)spaces endowed with  orthosymplectic forms.

\begin{theorem}
\label{th:main-1}
{\rm a)} For each
 contractive
 Lagrangian superlinear relation $P:V\arr \wt V$
there exists a linear
operator
$$
\B(P):\SF_{p,q}(\cA)\to\SF_{r,s}(\cA)
$$
such that

\smallskip

{\rm 1)} The following condition is satisfied
$$
\wh a(\wt h)\, \B(P)=\B(P)\,\wh a(h)
\q \text{for all $ h\oplus \wt h\in P$}.
$$

{\rm 2)} If $P$ is in the even component of Lagrangian super-Grassmannian, then
$\frB(P)$ is an integral operator with an even%
\footnote{see Subsection \ref{ss:antilinear-operators}.}
kernel. If $P$ is in the odd component, then $\B(P)\mathsf S$ is an
integral operator with an odd kernel.

\smallskip

Moreover, this operator is unique up to a scalar factor $\in \cA_{even}$.

\smallskip

{\rm b)} The operators $\B(P)$  are
 Gauss--Berezin operators
and all  Gauss--Berezin operators
arise in this way.
\end{theorem}

{\sc Proof.} Let us write 
  differential
equations for the kernel $K(z,\xi,\ov y,\ov\eta)$ of the operator
$\B(P)$.
  Denote
$$
h=v_+\oplus v_-\oplus w_+\oplus w_-,
\q
\wt h=\wt v_+\oplus \wt v_-\oplus \wt w_+\oplus \wt w_-
.
$$
Then
\begin{multline*}
\wh a( \wt h)
 \int K(z,\xi,\ov y,\ov\eta)\, f(y,\eta)\,
e^{-y\ov y^t-\eta\ov\eta^t}
 dy\,d\ov y\,d\ov\eta \,d\eta
 =\\=
  \int K(z,\xi,\ov y,\ov \eta)\,\wh a(h) f(y,\eta)\,
e^{-y\ov y^t-\eta\ov\eta^t}
 dy\,d\ov y\, d\ov \eta d\eta
.
\end{multline*}

Let $P$ be even.
Integrating by parts in the right-hand side, we get:
\begin{multline*}
\Bigl(\sum_i \wt v_+^{(i)}
\frac {\partial}{\partial z_i}+
\sum_i \wt v_-^{(i)} z_i
+
\sum_j \wt w_+^{(j)}
\frac {\partial}{\partial \xi_i}
+
\sum_j  \wt w_-^{(j)} \xi_j
\Bigr)K(z,\xi,\ov y,\ov\eta)
=\\=
\Bigl(\sum_i  v_+^{(i)}
\ov y_i+
\sum_i  v_-^{(i)} \frac {\partial}{\partial \ov y_i}
+
\sum_j  w_+^{(j)} \ov\eta_j
+
\sum_j   w_-^{(j)} \frac {\partial}{\partial\ov \eta_i}
\Bigr)K(z,\xi,\ov y,\ov \eta)
.
\end{multline*}
This system of equations
has the form (\ref{eq:41})
and determines a Gaussian.

\smallskip

Evenness condition was essentially used in this calculation. For
instance for an odd kernel $K$ we must write $(v_+^{(i)})^\sigma$
instead of $v_+^{(i)}$ in the right hand side.

\smallskip

Now, let $P$ be odd. Let us try to represent $\B(P)$
as a product
$$
\B(P)= \frC\cdot \frD[\eta_1] \cdot \mathsf S
.$$
Let $L$ be the kernel of $\frC$.
\begin{multline*}
\wh a( \wt h)
 \int L(z,\xi,\ov y,\ov\eta)\, \frD[\eta_1]\cdot \mathsf S\cdot f(y,\eta)\,
e^{-y\ov y^t-\eta\ov\eta^t}
 dy\,d\ov y\,d\ov\eta \,d\eta
 =\\=
  \int L(z,\xi,\ov y,\ov \eta)\,\frD[\eta_1]\cdot \mathsf S\cdot\wh a(h) f(y,\eta)\,
e^{-y\ov y^t-\eta\ov\eta^t}
 dy\,d\ov y\, d\ov \eta d\eta
.
\end{multline*}
Next, we change the order
$$
\frD[\eta_1]\, \mathsf\, S\wh a(v)
=
\wh a(\sigma(f_1) v)\, \frD[\eta_1]\, \mathsf S
,
$$
where $\sigma$ is an elementary reflection
of the type (\ref{eq:sigma-1})--(\ref{eq:sigma-2}).

We again get for $L$ a system of equations determining
a Gaussian.\hfill$\square$


\smallskip

{\bf\punct Products of Gauss--Berezin operators.%
\label{ss:main}}

\begin{theorem}
\label{th:main-2}
{\rm a)}
Let $P:\cA^{p|q}\arr \cA^{p'|q'}$, $Q:\cA^{p'|q'}\arr \cA^{p''|q''}$
be 
contractive Lagrangian relations. 
Assume that $P$, $Q$ are transversal. Then
\begin{equation}
\B(Q)\,\B(P)=\lambda\cdot \B(QP)
\label{eq:PQPQ}
,
\end{equation}
where $\lambda=\lambda(P,Q)$ is an even invertible phantom
constant.

\smallskip

{\rm b)} If $P$, $Q$ are not transversal, then
$$
\pia\bigl(\frB(Q)\,\frB(P) \bigr)=0
.
$$
\end{theorem}

{\sc Proof.} Let $v\oplus w\in P$, $w\oplus y\in Q$.
Then
$$
\B(Q)\B(P)\wh a(v)
=\B(Q)\wh a(w) \B(P)
= \wh a(y) \B(Q)\B(P)
.
$$
On the other hand,
$$
\wh a(y) \B(QP)=\B(QP)\wh a(v)
.
$$
By Theorem \ref{th:main-1},
 these relations define a unique
operator and we get
(\ref{eq:PQPQ}).

\smallskip

It remains to verify conditions of vanishing of
$$
\pia\bigl(\B(Q)\,\B(P)\bigr):\F_p\otimes \Lambda_q
\to\F_{p''}\otimes \Lambda_{q''}
.
$$
Here we refer to Subsection \ref{ss:remark}.
Our operator is a tensor product
of
\begin{equation}
\pia^+\bigl(\B(Q)\bigr)\,
\pia^+\bigl(\B(P)\bigr):\F_p
\to\F_{p''}
\label{eq:line1}
\end{equation}
and
\begin{equation}
\pia^-\bigl(\B(Q)\bigr)\,
\pia^-\bigl(\B(P)\bigr): \Lambda_q
\to \Lambda_{q''}
\label{eq:line2}
.
\end{equation}
The operator (\ref{eq:line1}) is
 a product of Gaussian integral operators.
By Theorem \ref{th:product-olshanski},
it is nonzero.

The operator (\ref{eq:line2}) is
 a product of Berezin operators.
It is nonzero if and only if $\pia^-(P)$ and $\pia^-(Q)$
are transversal, here we refer to Theorem
\ref{th:spinor-correspondence}.c.
\hfill $\square$

\begin{corollary}
For an element $g$ of  the Olshanski supersemigroup $\Gamma\OSp(2p|2q;\cA)$
{\rm(}see Subsect. {\rm\ref{ss:olshanski-super})} denote  its graph by $\mathrm{graph}(g)$.
Then $\frB(P(g))$ determines a projective representation
of  $\Gamma\OSp(2p|2q;\cA)$ over $\cA$,
$$
\frB(\mathrm{graph}(g_1))\,\frB(\mathrm{graph}(g_2))=\lambda(g_1,g_2)\frB(\mathrm{graph}(g_1g_2)), 
$$
where $\lambda(g_1,g_2)$ is an invertible element of $\cA$.
\end{corollary}

Denote by $G(2p|2q;\cA)$ the group of invertible elements
of Olshanski supersemigroup.
It is easy to see that $G(2p|2q;\cA)$
consists of all $g\in\OSp(2p|2q;\cA)$ such that $\pia(g)\in\Sp(2q,\R)\times \OO(2p,\C)$%
\footnote{In $\Sp(2q,\R)\times \OO(2p,\C)$, we have a product of a real Lie group and a  complex Lie group,
so it is not a real form (on real forms, see \cite{Ser}, \cite{FG}).
 Moreover, $G(2p|2q;\cA)$ is not a supergroup
in the sense of the usual definitions \cite{BK}, \cite{DeW}, \cite{CF}, \cite{Man}, \cite{BL} 
(because there are no intermediate 
Lie superalgebras between $\fro\frs\frp(2p|2q,\R)$ and $\fro\frs\frp(2p|2q,\C)$). However, 
$G(2p|2q;\cA')$ depends functorially on algebras $\cA'$ described in Subsection \ref{ss:super}.
}.

\begin{corollary} The map
$g\mapsto\frB(\mathrm{graph}(g_1))$ determines a projective representation of the group
$G(2p|2ql\cA)$ over $\cA$.
\end{corollary}


\section{Final remarks}

\label{s:howe}

\COUNTERS


\smallskip

{\bf\punct Extension of notion of Gaussian operators?%
\label{ss:?}}
Our main result seems incomplete, since
Theorem \ref{th:main-2}
describes product of integral
operators only if superlinear relations are
transversal. However, a product of integral operators
  can be written explicitly in all  cases
 (see Subsection \ref{ss:second-way}).
We arrive at the following question:

\smallskip

\begin{question}
 {\rm a)} Is it possible to extend
the definitions of Gaussian operators
and Lagrangian superlinear relations to make the
formula \eqref{eq:PQPQ} valid for all $P$, $Q${\rm ?}

\smallskip

{\rm b)} Is it reasonable to consider the expressions
\eqref{eq:general} as Gaussians{\rm ?}
\end{question}


\smallskip

{\bf\punct The infinite-dimensional orthosymplectic supergroup.}
Recall that constructions of orthogonal and symplectic spinors described 
in Sect. \ref{s:orthogonal}--\ref{s:symplectic}  survive well in the infinite-dimensional limit,  \cite{Ner-book}, Chapters IV, VI. However,
there remain gaps between
sufficient and necessary conditions for boundedness of Gaussian
operators  in the bosonic Fock spaces 
with infinite number of degrees of freedom
(see  \cite{Ner-book}, Sect. VI.3-4, \cite{Olsh-gauss}) and
similar gaps in fermionic case (see \cite{Ner-book}, Sect. IV.2).
On the other hand, even for finite-dimensional supergroups 'unitary
representations' are realized by unbounded operators in super-Hilbert spaces
(for instance, for orthosymplectic spinors discussed in this work).

\begin{question}
Find natural topologies in Fock--Berezin space 
$\F_{\infty,\infty}(\cA)$ with infinite number of variables $z_1$, $z_2$, \dots, 
$\xi_1$, $\xi_2$, \dots and conditions of boundedness of Gauss--Berezin operators.
\end{question}

Such topologies can depend on further applications. See, for instance, the next subsection, a straightforward repetition of our definition of $\boldsymbol{\cS}\F_{p,q}$ does not work in that situation.

\sm

{\bf\punct The super-Virasoro algebras.}
It is well-known that highest weight representations of the Virasoro
algebra appear in a natural way as restrictions of infinite-dimensional Weil representation and infinite-dimensional spinor representation to the Lie algebra of vector fields on the circle, see \cite{Seg}, \cite{Ner-highest}, for more details, see \cite{Ner-book}, Chapter VII. Let us explain that a similar phenomenon takes place for super-Virasoro algebras
and ortho-symplectic spinors. 

 Consider the circle
$S^1$ with the coordinate $\phi\in \R/2\pi \Z$, denote by $C_+^\infty(S^1)$
(resp., $C_-^\infty(S^1)$) the space
of smooth complex functions on $S^1$ such that $f(\phi+\pi)=f(\phi)$ (resp.,  $f(\phi+\pi)=-f(\phi)$).
Let us denote by $f(\phi)(d\phi)^\lambda$  densities of weight $\lambda$ on the circle,
$d\phi=  (d\phi)^1$, $\frac d{d\phi}=(d\phi)^{-1}$.
Recall that vector fields act on densities by 
$$
  a(\phi)\frac d{d\phi} \Bigl(f(\phi)(d\phi)^{\lambda} \Bigr)
  =\Bigl(a(\phi) f'(\phi)+\lambda a'(\phi) f(\phi) \Bigr) (d\phi)^{\lambda}.
  $$
  Then, $d\phi=(d\phi)^1$, $\partial/\partial \phi=(d\phi)^{-1}$.

Consider a 'super-Witt' Lie superalgebra $\mathfrak{sw}_-$ (see Kirillov \cite{Kir}), whose elements are
\begin{equation}
a(\phi) \frac d{d\phi}\,\oplus \, b(\phi) (d\phi)^{-1/2},
\qquad a\in C_+^\infty(S^1),\, b\in C_-^\infty(S^1),
\label{eq:a-b}
\end{equation}
  the first summand has parity $\ov 0$, the second summand parity
  $\ov 1$. The supercommutator is defined as follows: 
  
  \sm
  
 --- supercommutator of vector fields is the usual commutator;
 
 \sm
 
 --- supercommutator of a vector field and a density is the natural action of  vector fields  on densities
 of weight $-1/2$;

  --- the anticommutator
  of densities $b_1(\phi)(d\phi)^{-1/2}$, $b_2(\phi)(d\phi)^{-1/2}$
  is the product $2b_1(\phi)b_2(\phi) \frac d{d\phi}$.
  
  \sm
  
  This supercommutator is a 'natural
  differential geometric operation%
  \footnote{On natural differential operations, see, e.~g.,  \cite{Kol}, \cite{Gro}, \cite{GLS}.},
  i.e., it is invariant with respect to the action of the group of diffeomorphisms $r$ of $S^1$ satisfying $r(\phi+\pi)=r(\phi)+\pi$.
  
 Consider the superlinear space%
 \footnote{Cf. constructions for Virasoro algebra and the group of diffeomorphisms of circle in \cite{Ner-book}, Subsect.VII.2.2-2.3, 2.5, VII.3, Sect. VII.3, VIII.6.}
  $W:=W_{\ov 0}\oplus W_{\ov 1}$ consisting of 
 \begin{equation}
 f(\phi)\,\oplus\,  g(\phi) (d\phi)^{1/2}, 
\qquad f\in C_+^\infty(S^1)/\C,\, g\in C_-^\infty(S^1).
\label{eq:f-g}
 \end{equation}
 We equip $W$  with the orthosymplectic form
  by
 \begin{multline*}
 \bigl\{ f_1(\phi)\,\oplus\,  g_1(\phi) (d\phi)^{1/2},\,
  f_2(\phi)\,\oplus\,  g_2(\phi) (d\phi)^{1/2} \bigr\}=\\=
 \frac 14\int_{S^1} (f_1 df_2-f_2 df_1)+ \int_{S^1} g_1 g_2\, d\phi.
 \end{multline*}
 Notice, that the element $1\oplus 0$ is contained in the kernel of this form and we really get a form on the quotient $W:=(C_+^\infty(S^1)/\C)\oplus C_-^\infty(S^1)$.
 Next, we define the inner product on $W$ by
 \begin{multline*}
\bigl\la \sum p_n e^{2in\phi}\oplus\sum q_n  e^{(2n+1)\phi} (d\phi)^{1/2},
\sum p_n' e^{2in\phi}\oplus\sum q_n'  e^{(2n+1)\phi} (d\phi)^{1/2}
\bigr\ra:=\\=
\sum |n| \, p_n \ov{p_n'}+\sum q_n \ov{q_n'}.
 \end{multline*}
The superalgebra $\mathfrak{sw}_-$ acts in the space $W$ as follows.
Vector fields act in $W_{\ov 0}\oplus W_{\ov 1}$ in a natural way.
$$
a(\phi)\frac{d}{d\phi}\Bigl(f(\phi)\oplus g(\phi)d\phi^{1/2}\Bigr)= a(\phi)f'(\phi)\oplus \Bigl(a(\phi) b'(\phi)+\tfrac12 a'(\phi) b(\phi) \Bigr) (d\phi)^{1/2}.
$$
An element $0\oplus b(\phi)(d\phi)^{-1/2}$
acts as 
\begin{multline*}
f(\phi)\oplus g(\phi)\mapsto b(\phi)g(\phi)\,\oplus\, b(\phi) f'(\phi)(d\phi)^{1/2}=\\=
b(\phi)(d\phi)^{-1/2}\cdot g(\phi)(d\phi)^{1/2}\,\oplus\,  df(\phi)\cdot b(\phi)(d\phi)^{-1/2} .
\end{multline*}
Again, this action is a 'natural differential geometric operation'.

In this way, we get   an embedding of $\mathfrak{sw}_-$ to an infinite-dimensional Lie  superalgebra $\mathfrak{osp}(2\infty|2\infty)$. 
Next, we apply the infinitesimal version of the orthosymplectic spinors setting $p=\infty$, $q=\infty$ 
 in the notation Subsect. \ref{ss:1.1} and assuming even variables $x_j$ be complex.
 We arrive at a {\it projective} representation of the Lie superalgebra $\mathfrak{sw}_-$
 in a tensor product of the bosonic and fermionic Fock spaces with infinite number degrees of freedom. Its restriction to the Lie superalgebra $\mathfrak{sw}_-$ 
also is a projective representation. The corresponding central extension 
 is 
 the {\it Neveu--Schwarz super-Virasoro algebra}. For formulas for operators,  see, e.g., \cite{Gre}, Subsect. 4.2.2.
 
 \sm
 
{\sc Remark.} 
 The {\it Ramond super-Virasoro algebra} arises in a similar way: in \eqref{eq:a-b} we assume
 $a$, $b\in C^\infty(S^1)$  and get another 'super-Witt' algebra 
 $\mathfrak{sw}_+$.
 Next, in \eqref{eq:f-g} we assume  $f$, $g\in C^\infty(S^1)$. 
 Repeating the construction above we get an embedding of $\mathfrak{sw}_+$ to an infinite-dimensional Lie superalgebra $\mathfrak{osp}(2\infty|2\infty+1)$, a construction of the spinor representation
 for this algebra must be slightly modified (cf. \cite{Ner-book}, Subsect. III.3.4, VII.3.4). We get a projective representation of $\mathfrak{sw}_+$
 or a linear representation of the Ramond Lie superalgebra.
 \hfill $\boxtimes$

 \begin{question}
 {\rm a)} Integrate these representations of the Neveu--Schwarz and Ramond 
 algebras to  actions of the corresponding supergroups
 by Gauss--Berezin operators.
 
 \sm
 
 {\rm b)} Extend
highest weight representations of the Ramond and Neveu--Schwarz
 supergroups to  complex semigroups as in 
{\rm \cite{Ner-holom}, \cite{Ner-book}, Sect. VII.4-5}. 
 \end{question}
 
 In this context, it is more natural to consider a more general family of super-Virasoro algebras, see 
 \cite{Kac}, \cite{GLS}.


\smallskip

{\bf\punct  'Unitary representations' of supergroups.%
\label{ss:??}}
  The  Howe duality for orthosymplectic spinors exists and 
  it was a subject of numerous works,
e.g.,
 \cite{Nish}, \cite{China}, \cite{China2}, \cite{LI},
\cite{Cou}.
In particular, this automatically produces  many representations
of supergroups, which can be extended  to Grassmannians.

As far as I know, a notion of {\it a unitary representation of a Lie supergroup} is commonly recognized. Consider a finite-dimensional {\it real} Lie superalgebra 
$\frg=\frg_{\ov 0}\oplus \frg_{\ov 1}$, let $G_{\ov 0}$ be the Lie group, corresponding
to $\frg_{\ov 0}$. We consider a super-Hilbert space
 $H=H_{\ov 0}\oplus H_{\ov 1}$ and a representation of $\frg$ such that 
 the action of $\frg_{\ov 0}$ corresponds to a unitary representation of $G_{\ov 0}$, and for each $X\in \frg_{\ov 1}$ the operator $\sqrt[4]{-1}X$ is self-adjoint.

 Under these conditions, for each  $X\in \frg_{\ov 1}=(\sqrt[4]{-1}X)^2$
  the operator 
 $\sqrt{-1} X^2$  is positive (semi)definite.
On the other hand, $X^2=\frac 12[X,X]_{s}$ is a generator of $\frg_{\ov 0}$. 
  In unitary representations, generators of Lie algebras rarely have spectra supported by a semi-axis. 
 This places strong constrains on the Lie algebra $\frg_{\ov 0}$ and its representation $\tau$ in $H$. In
 the semi-simple case this implies (see \cite{Olsh-semi}) that $\tau$
 is an element of Harish-Chandra \cite{Harish} highest weight holomorphic series. 
 In, particular, few semisimple real supergroups have non-trivial unitary representations in this sense%
 \footnote{Both highest weight and lowest weight representations of a supergroup, say $\OSp(2p|2q;\R)$, can be unitary, but they correspond to different $\sqrt[4]{-1}$ (see \cite{Fur1}). For this reason, 
 a tensor product of a highest weight and a lowest weight unitary representation
 is not unitary.
Counterparts of principal series representations for supergroups exist (since there are flag supervarieties and line bundles over them) but they are not unitary. This restricts possible ways for applications of Lie superalgebras to special
functions. 
This is also an obstacle for an extension
of Olshanski theory \cite{Olsh-GB} of representations of infinite-dimensional real classical groups to supergroups. 
\newline
So a picture is poor comparatively to classical representation theory.
It remains a minor hope that a notion of  unitary representations 
of supergroups admits some extension.}.
 This  was observed at least in works by Furutsu and Nishiyama \cite{Fur1}, \cite{Fur2}.
 Such 'unitary representations' of superalgebras are nice counterparts  of the Harish-Chandra holomorphic series 
 and admit explicit classification (Jakobsen \cite{Jak}), which is a counterpart of the Howe--Enright--Wallach classification
 of unitary highest weight representations.

\begin{conjecture}
Any unitary representation  of a classical Lie supergroup admits an extension
to a certain domain in a certain Grassmanian as in our Theorems
{\rm \ref{th:main-1}--\ref{th:main-2}}.
\end{conjecture}

\noindent
{\tt Math.Dept., University of Vienna;\\
Institute for Theoretical and Experimental Physics\\
 (until 11.2021);\\
Insitute for Information Transmission Problems;\\
 Moscow State University, MechMath Department.\\
URL:www.mat.univie.ac.at/$\sim$neretin

}

\end{document}